\documentclass{amsart}
\ifdefined\pdfinfoomitdate\pdfinfoomitdate=1\fi
\ifdefined\pdfsuppressptexinfo\pdfsuppressptexinfo=-1\fi
\ifdefined\pdftrailerid\pdftrailerid{}\fi

\usepackage{amssymb,mathtools}
\usepackage{multirow}
\usepackage{algorithm}
\usepackage[noend]{algpseudocode}
\usepackage{placeins,needspace}
\usepackage{aliascnt}
\usepackage{tikz}
\usetikzlibrary{arrows.meta,calc}
\definecolor{omniBlue}{RGB}{37,86,126}
\definecolor{omniOrange}{RGB}{180,95,35}
\usepackage[hidelinks]{hyperref}
\hypersetup{hidelinks,pdfauthor={},pdftitle={},pdfsubject={},pdfkeywords={},
  pdfcreator={},pdfproducer={},pdfcreationdate={},pdfmoddate={},pdftrapped={}}
\usepackage[initials]{amsrefs}
\renewcommand{\PrintDOI}[1]{DOI~\href{https://doi.org/#1}{\nolinkurl{#1}}}
\usepackage[capitalize,nameinlink]{cleveref}

\makeatletter
\renewcommand{\theALG@line}{\thealgorithm.\arabic{ALG@line}}
\renewcommand{\theHALG@line}{\thealgorithm.\arabic{ALG@line}}
\makeatother

\numberwithin{equation}{section}

\theoremstyle{plain}
\newtheorem{theorem}{Theorem}[section]
\newaliascnt{lemma}{theorem}
\newtheorem{lemma}[lemma]{Lemma}
\aliascntresetthe{lemma}
\newaliascnt{proposition}{theorem}
\newtheorem{proposition}[proposition]{Proposition}
\aliascntresetthe{proposition}
\newaliascnt{corollary}{theorem}
\newtheorem{corollary}[corollary]{Corollary}
\aliascntresetthe{corollary}
\theoremstyle{definition}
\newaliascnt{definition}{theorem}
\newtheorem{definition}[definition]{Definition}
\aliascntresetthe{definition}
\newaliascnt{assumption}{theorem}

\aliascntresetthe{assumption}
\theoremstyle{remark}
\newaliascnt{remark}{theorem}

\newtheorem*{remark*}{Remark}
\aliascntresetthe{remark}

\crefname{theorem}{Theorem}{Theorems}
\Crefname{theorem}{Theorem}{Theorems}
\crefname{lemma}{Lemma}{Lemmas}
\Crefname{lemma}{Lemma}{Lemmas}
\crefname{proposition}{Proposition}{Propositions}
\Crefname{proposition}{Proposition}{Propositions}
\crefname{corollary}{Corollary}{Corollaries}
\Crefname{corollary}{Corollary}{Corollaries}
\crefname{definition}{Definition}{Definitions}
\Crefname{definition}{Definition}{Definitions}
\crefname{remark}{Remark}{Remarks}
\Crefname{remark}{Remark}{Remarks}

\newcommand{\N}{\mathbb N}
\newcommand{\Nzero}{\mathbb N_0}
\newcommand{\R}{\mathbb R}
\newcommand{\1}{\mathbf 1}
\newcommand{\Kcal}{\mathcal K}
\newcommand{\Tcal}{\mathcal T}
\newcommand{\Fcal}{\mathcal F}
\newcommand{\Hcal}{\mathcal H}
\newcommand{\Pcal}{\mathcal P}
\newcommand{\Lcal}{\mathcal L}
\newcommand{\Ical}{\mathcal I}
\newcommand{\Mcal}{\mathcal M}
\newcommand{\Acal}{\mathcal A}
\newcommand{\rootkey}{\kappa_0}
\newcommand{\leaf}{\mathsf{leaf}}
\newcommand{\ch}{\operatorname{ch}}
\newcommand{\parop}{\operatorname{par}}
\newcommand{\nodes}{\operatorname{nodes}}
\newcommand{\leaves}{\operatorname{leaves}}
\newcommand{\tr}{\operatorname{tr}}
\newcommand{\relint}{\operatorname{relint}}
\newcommand{\card}[1]{\left\lvert #1\right\rvert}
\newcommand{\set}[1]{\left\{#1\right\}}
\newcommand{\given}{\,:\,}
\newcommand{\dotcup}{\mathbin{\dot\cup}}
\newcommand{\Reconstruct}{\mathsf N}
\newcommand{\Apply}{\mathsf A}
\newcommand{\anc}{\operatorname{anc}}

\begin{document}

\title[Forests of Omnitrees]{Forests of Omnitrees: Analysis and Applications to Anisotropic Space--Time Finite Elements}
\author{Nils Margenberg}
\address{Institute of Analysis and Numerics, Otto von Guericke University Magdeburg,
Universitätsplatz 2, 39106 Magdeburg, Germany}
\email{nils.margenberg@ovgu.de}
\date{}

\begin{abstract}
We prove output-optimal work bounds for omnitree forests, which bisect
cells along selected axes. A canonical descriptor characterizes representable
partitions independently of split order. Complete one-level cycles, including
coarsening, face balance and selected coarse meshes, take sequential work
proportional to stored input and output at fixed dimension and key width,
with polynomial depth dependence even for unbalanced input. From balanced
input, refinement-only balance stays within each request's coordinate slab;
this bound is attained. Polynomial trace equations characterize conforming
high-order spaces across crossed face partitions. Reconstruction moments
bound interpolation errors and capture the logarithmic loss that bilinear
nodal interpolation can suffer on balanced meshes. Elements of degree up to
four resolve moving heat layers on unstructured 4D space-time and 5D
space-time-parameter meshes. Multigrid uses the generated hierarchies.
\end{abstract}

\maketitle

\section{Introduction}
\label{sec:introduction}
In space-time and parametric discretizations, a moving front becomes a thin
layer that needs fine resolution across its thickness and less resolution along
it.
Anisotropic refinement bisects cells along selected axes
\cites{neumuellerSteinbachFlexibleSpaceTime2011,
apelAnisotropicFiniteElements1999,bungartzGriebelSparseGrids2004}:
a one-axis bisection creates two children. An all-axis split creates 16
in 4D and 32 in 5D.
Omnitree forests enable local refinement of unstructured space-time and
space-time-parameter meshes with conforming high-order finite elements,
tensor-product kernels and automatically generated coarse spaces for multigrid.
PDEs with random coefficients \cite{daworStochasticGalerkin2026},
as well as inverse and optimal control problems,
motivate resolving space, time and parameter dependence together.

The computational problem is to update these meshes, find neighboring and
overlapping cells, and build coarse multigrid meshes with work proportional
to explicitly stored input and output.

An \emph{omnitree} bisects a coarse \emph{macro cell} along selected axes,
extending quadtrees and octrees
\cites{sametQuadtreeRelatedHierarchical1984,bursteddeP4estScalableAlgorithms2011}.
Different split orders can produce the same cells. Geometrically
mergeable cells need not share a parent. Merging them can leave no tree
representation. A canonical descriptor characterizes representable partitions
independently of split order (\Cref{thm:normal-form}).
Neighboring face partitions can cross, like vertical and horizontal strips,
leaving neither side as a coarse
interpolation trace.

In each coordinate direction, an adaptation step permits at most one level
of requested refinement or coarsening, with refinement taking precedence.
Each cycle applies requested refinement and permitted coarsening,
balances tangential face resolutions and optionally removes crossed traces.
It emits the final forest, nested coarse
meshes and \emph{cell incidence}: positive-volume overlaps between meshes.

\begin{figure}[!hbp]
\centering
\begin{tikzpicture}[x=1cm,y=1cm,>=Stealth,
    every node/.style={font=\small},
    stage/.style={draw=black!55,rounded corners=1pt,line width=.6pt,
      text width=2.45cm,minimum height=.8cm,inner sep=3pt,align=center},
    cycle/.style={stage,minimum height=1.05cm}]
  \node[anchor=west,font=\small\bfseries] at (-1.35,.75) {Reading order};
  \node[anchor=east,font=\footnotesize] at (10.8,.75)
    {Tree constructions, then face interactions and accounting};
  \node[stage,fill=black!2] (representation) at (0,0)
    {\Cref{sec:foundations,sec:normal-form}\\Representation};
  \node[stage,draw=omniBlue!85!black,fill=omniBlue!4] (adaptation) at (3.15,0)
    {\Cref{sec:normalization-adaptation}\\Adaptation};
  \node[stage,fill=black!2] (ranks) at (6.3,0)
    {\Cref{sec:hierarchy}\\Rank cuts};
  \node[stage,draw=omniOrange!85!black,fill=omniOrange!4] (faces) at (9.45,0)
    {\Cref{sec:faces}\\Face closures};
  \draw[->,line width=.7pt] (representation.east)--(adaptation.west);
  \draw[->,line width=.7pt] (adaptation.east)--(ranks.west);
  \draw[->,line width=.7pt] (ranks.east)--(faces.west);
  \node[anchor=west,font=\small\bfseries] at (-1.35,-.7)
    {Mesh-cycle order (\Cref{sec:composed-accounting}): reuse these constructions};
  \node[cycle,draw=omniBlue!85!black,fill=omniBlue!4] (adapt) at (0,-1.55)
    {Directional adaptation\\
     {\footnotesize \Cref{sec:normalization-adaptation}}};
  \node[cycle,draw=omniOrange!85!black,fill=omniOrange!4] (closure) at (3.15,-1.55)
    {Face closures\\
     {\footnotesize \Cref{sec:faces}}};
  \node[cycle,fill=black!2] (emit) at (6.3,-1.55)
    {Canonical reconstruction\\
     {\footnotesize \Cref{sec:normal-form}}};
  \node[cycle,fill=black!2] (hierarchy) at (9.45,-1.55)
    {Selected rank cuts\\
     {\footnotesize \Cref{sec:hierarchy}}};
  \draw[->,line width=.7pt] (adapt.east)--(closure.west);
  \draw[->,line width=.7pt] (closure.east)--(emit.west);
  \draw[->,line width=.7pt] (emit.east)--(hierarchy.west);
\end{tikzpicture}
\par\vspace{6pt}
\caption{Reading order (top) and mesh-cycle order (bottom).
Sections~\ref{sec:foundations}--\ref{sec:hierarchy} develop tree constructions
independently of face closure. \Cref{sec:faces} combines them with face
interactions to prove \Cref{cor:flat-cycle-accounting}.}
\label{fig:cycle-work-map}
\end{figure}
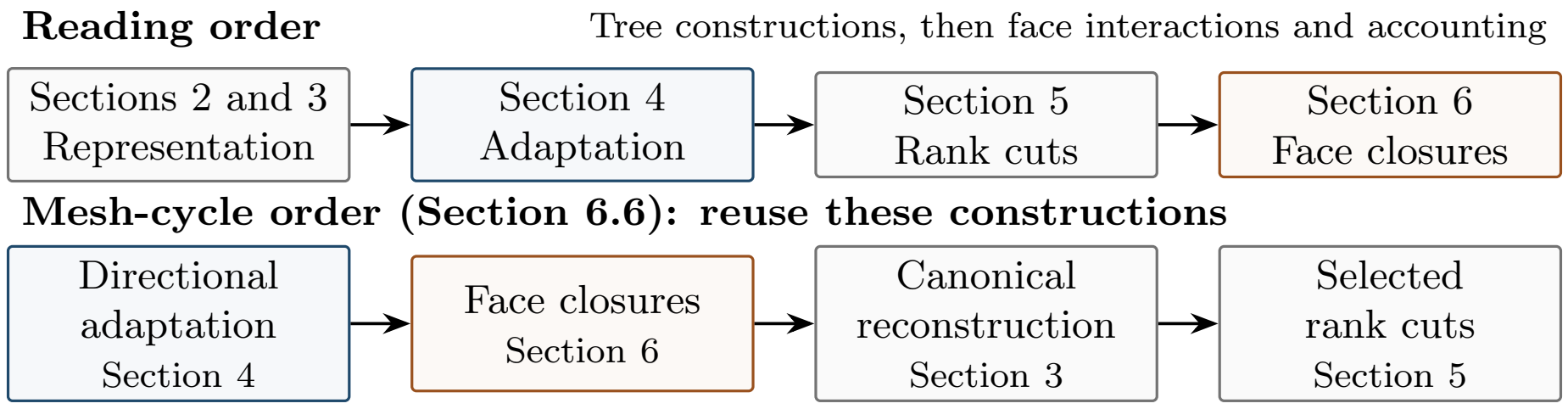

The main work bound covers this complete mesh cycle (\Cref{fig:cycle-work-map}). Here $n,q,n_*$ count
input descriptor records, requests and final descriptor records, and $N_j$
counts records at each selected hierarchy level $j\in J$.
At fixed dimension $d$ and key width $w\ge L_{\max}$,
the mesh cycle is output-optimal in flat sequential work
(\Cref{cor:flat-cycle-accounting}):
\[
 T_{\rm cycle}=O_{d,w}\!\Big(n+q+n_*+\sum_{j\in J}N_j\Big).
\]
The bound includes repeated split reordering and contact updates. Its explicit
form \eqref{eq:cycle-gap-depth-cost} has polynomial depth dependence even
for unbalanced input.

Balance can turn one requested split into many. For refinement-only cycles
from balanced input, propagation follows a \emph{sheet}: the request's
coordinate slab carried across macro faces. Only cells spanning that slab in
the transported direction are reached. A layer family attains this geometric
bound. Several requests can require the same cell to split in the same
direction; the growth bound counts this once. Depth weights are linear in 2D
and necessarily exponential in higher dimensions for this whole-cell rule
(\Cref{thm:marked-growth}).

Polynomial trace equations characterize conforming spaces across crossed
partitions (\Cref{thm:conforming-kernel}). Reconstruction moments bound
interpolation errors by weighted distances to the samples that determine
constrained nodal values. Even on balanced meshes with nested face
partitions, standard bilinear nodal interpolation can lose a logarithmic
factor in the local $h^2$ estimate (\Cref{prop:local-interpolation-loss}).

The term ``omni-tree'' appears in the SPLITFLOW CFD report
\cite{jordanSplitflow1998}*{Sec.~2}. Earlier dyadic CART
\cite{donohoCart1997}*{Sec.~2.1} generates the same two-dimensional leaf
partitions: successive one-axis splits replace simultaneous cuts.
Directional trees also appear in Cartesian PDE solvers
\cite{ogawaAnisotropicCartesian2003} and high-dimensional adaptive meshes
\cite{grandinHighDimensional2015}.

Equivalent rectangular patterns can have different cut sequences
\cite{brinHigherDimensional2004}*{Sec.~2.1}. Common-axis lifting is established in
\cite{hennigMatucciPresentations2012}*{Lemma~6, Corollary~9, Proposition~10}.
Normalized omnitrees and refinement estimates appear in
\cite{pollingerBeautyAnisotropic2026}, with coarsening and restructuring in
\cite{pollingerFullyDynamic2026}.
Our work bounds cover canonical regrouping, admissible directional reduction,
nested rank cuts and face closure as a complete cycle.

Forest meshes separate macro connectivity from refinement. p4est supplies
scalable isotropic refinement \cite{bursteddeP4estScalableAlgorithms2011},
balance \cite{isaacBursteddeGhattasBalance2012} and recursive traversal with
FE node construction \cite{isaacRecursiveOctrees2015}.
Algorithms for hybrid forests are developed in
\cite{holkeGhostHybrid2021}.
Grid abstractions \cite{bastianGridInterface2008} and adaptive FE infrastructure
\cite{bangerthGenericAdaptive2011} connect refinement to discretization.
p6est supports anisotropic ice-sheet meshes \cite{isaacStadlerGhattasStokes2015}.
Automatic anisotropic 3D $hp$ refinement is studied in
\cite{rachowiczPardoDemkowiczHp2006}.
Nonconforming 2D/3D methods support anisotropic refinement, derefinement and
unlimited neighboring refinement ratios. Their high-order interpolation uses
further refinement to establish master/slave relations on crossed traces
\cite{cervenyNonconforming2019}*{Sec.~5.3}. Our two-sided equations retain
crossed traces on active meshes and hierarchy cuts, making trace-comparability
closure optional.

The hyper.deal framework forms tensor products of lower-dimensional meshes for matrix-free
methods in up to six dimensions \cite{munchHyperdeal2021}. Omnitree meshes
retain tensor-product bases without requiring a spatial, temporal or parameter
factorization. Their trace equations also allow directional polynomial degrees.

Space-time parabolic formulations include
\cites{steinbachSpaceTimeFinite2015,langerSpaceTime2015,langerSchafelnerAdaptive2020};
anisotropic hexahedral estimates and hanging-node refinement appear in
\cite{langerSchafelnerHexahedral2022}. Prismatic refinement admits cumulative
completion bounds and stable quasi-interpolation with prescribed space-time
scaling and controlled hanging nodes
\cite{morinSchneiderSchneiderApproximation2026}; omnitrees permit independent
coordinate refinement. Related solvers include space-time multigrid
\cite{ganderNeumuellerSpaceTime2016}, hyperbolic discontinuous Galerkin methods
\cite{dorflerFindeisenWieners2016}, and temporal wavelets with adaptive low-rank
spatial approximations \cite{bachmayrFaldumSpaceTime2024}.
Hierarchical bases and $hp$-multigrid support adaptation
\cites{mitchellHierarchicalRefinement1991,mitchellHpMultigrid2010}.
Our rank cuts supply coarse meshes on directional forests.

\Cref{sec:foundations} fixes the mesh and work model.
\Cref{sec:normal-form,sec:normalization-adaptation,sec:hierarchy,sec:faces}
develop representation, adaptation, hierarchies, face closure and their work bounds.
\Cref{sec:finite-element-realization,sec:numerical-experiments,sec:conclusion}
give the FE realization, experiments and conclusions.
Appendix~\ref{sec:implementation} contains auxiliary proofs.

\section{Mathematical setting}
\label{sec:foundations}
\label{sec:keys}

\subsection{Dyadic cells}

Integer keys record a cell's reference position and width in each axis,
independently of tree grouping.
We work in dimension \(d\ge1\), write \([d]=\{1,\ldots,d\}\), and use an arbitrary finite
coordinate-level bound \(L_{\max}\in\N\).
Write $I_{\ell,i}=[i2^{-\ell},(i+1)2^{-\ell})$.
The finite key set and its product cells are
\begin{equation}
\label{eq:dyadic-key-cell}
 \begin{gathered}
 \Kcal_d:=\{\kappa=((\ell_a,i_a))_{a=1}^d\in\Nzero^{2d}:
       0\le\ell_a\le L_{\max},\ 0\le i_a<2^{\ell_a}\},\\
 Q_\kappa:=\prod_{a=1}^d I_{\ell_a,i_a}\subset[0,1)^d.
 \end{gathered}
\end{equation}
Write $\ell(\kappa)=(\ell_a)_a$ for the level vector.
For example, $((1,0),(2,3))$ identifies $[0,\tfrac12)\times[\tfrac34,1)$.
All requested splits remain in \(\Kcal_d\). The root key is
$\rootkey=((0,0),\ldots,(0,0))$.
The dyadic ancestor order is \(\kappa\preceq\eta\) if \(Q_{\eta}\subseteq Q_\kappa\). Equivalently,
\begin{equation}
\label{eq:containment}
  \ell_a(\kappa)\le \ell_a(\eta),
  \qquad
  i_a(\kappa)=\left\lfloor 2^{-(\ell_a(\eta)-\ell_a(\kappa))}i_a(\eta)\right\rfloor
  \quad (a\in[d]).
\end{equation}

Two dyadic intervals are disjoint or one contains the other. Hence a
nonempty product-cell intersection takes the smaller interval in each coordinate:
\begin{equation}
\label{eq:dyadic-dichotomy}
 Q_\kappa\cap Q_\eta
 =\prod_{a=1}^d\bigl(I_{\ell_a(\kappa),i_a(\kappa)}
                  \cap I_{\ell_a(\eta),i_a(\eta)}\bigr)=Q_\xi.
\end{equation}
The cells themselves need not be nested: one can be narrower in $x$,
the other in $y$.

A split mask is a set of axes to bisect: $\sigma\subseteq[d]$.
For nonempty \(\sigma\) and branch bits \(b\in\{0,1\}^{\sigma}\), define the child
key \(c_{\sigma,b}(\kappa)\) by
\begin{equation}
\label{eq:child}
  (\ell_a',i_a')=
  \begin{cases}
    (\ell_a+1,2i_a+b_a),& a\in\sigma,\\
    (\ell_a,i_a),& a\notin\sigma .
  \end{cases}
\end{equation}
Write $\ch_\sigma(\kappa)=\{c_{\sigma,b}(\kappa):b\in\{0,1\}^{\sigma}\}$.
By \eqref{eq:child}, these children partition the parent:
\begin{equation}
\label{eq:child-partition}
  Q_\kappa=\dotcup_{\eta\in\ch_\sigma(\kappa)} Q_\eta,
  \qquad
  \card{\ch_\sigma(\kappa)}=2^{\card\sigma} .
\end{equation}

Containment, equality, overlap, and child construction cost $O(d)$
integer operations.

\subsection{Omnitrees and forests}
\label{sec:trees}

A tree partitions a reference cell. A forest is a collection of these trees
over a macro complex (\Cref{sec:macro-complexes}).

\subsubsection{Trees and descriptors}

An omnitree is a recursive partition tree: each node is a leaf or has children
from a split along selected axes. A descriptor records this recursion:
\[
  U=\leaf,
  \qquad\text{or}\qquad
  U=(\sigma,(U_b)_{b\in\{0,1\}^{\sigma}}),
  \quad \sigma\ne\emptyset,
\]
where $U_b$ describes child $b$. The root mask is
$\operatorname{rm}(\leaf)=\emptyset$ and
$\operatorname{rm}(U)=\sigma$ for an internal node.
\emph{Preorder} visits a node before its children, recursively in branch
order. \emph{Postorder} visits the children first. Preorder stores the masks
as one sequence \cite{pollingerBeautyAnisotropic2026}*{Section~4.1}.

\begin{definition}[Omnitree descriptor]
\label{def:descriptor}
A stored descriptor is a finite sequence of split masks in preorder:
$\emptyset$ for a leaf, or a nonempty mask followed by its child descriptors
in increasing binary order, with the lowest split axis as the least
significant bit. It is valid at a root key if recursive reading consumes
the sequence exactly and every induced key lies in $\Kcal_d$.
\end{definition}

Because masks fix child counts, the sequence $U$ decodes uniquely into a tree
$\Tcal$, with one descriptor record for each node, internal or leaf.
The relation $v\mapsto\kappa(v)\mapsto Q_{\kappa(v)}$ connects each node to
its cell. Leaf nodes identify the active cells.
The simultaneous split in \Cref{fig:normal-form-commutation} is stored as
$(\{x,y\},\emptyset,\emptyset,\emptyset,\emptyset)$: five records, four leaves.
Write $\ch(v)$ and $\parop(v)$ for node children and parent.
$\ch_\sigma(\kappa)$ in \eqref{eq:child-partition} instead contains child keys.

\subsubsection{Partitions}

Tree leaves are the active cells. A partition is \emph{admissible} if a
descriptor produces it using \eqref{eq:child-partition}. Set
$\Pcal(\Tcal):=\{Q_{\kappa(v)}:v\in\leaves(\Tcal)\}$.

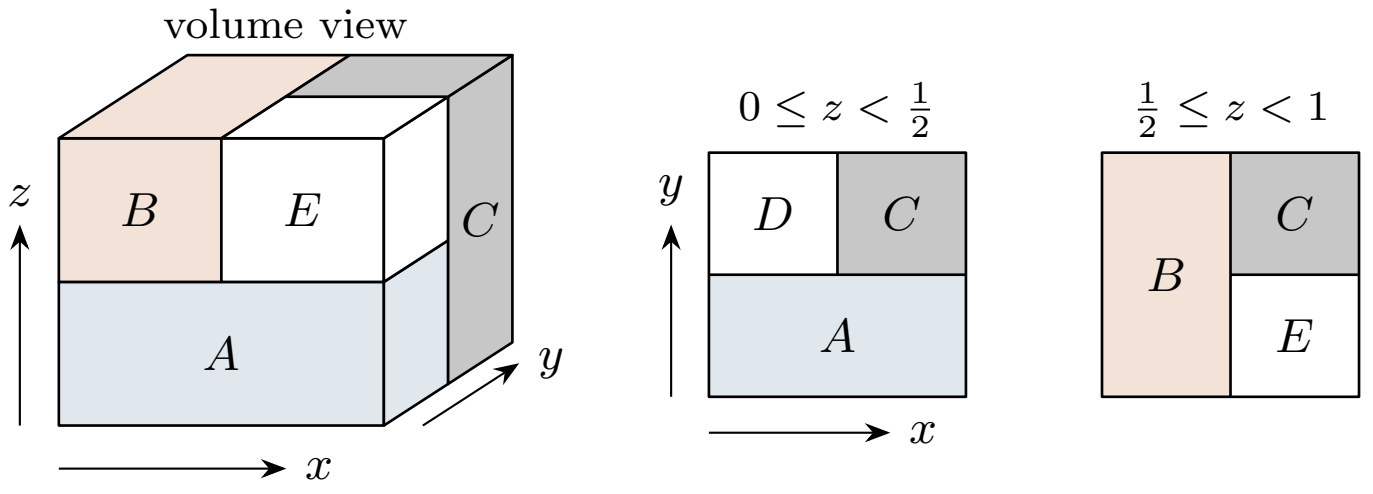
\begin{figure}[htbp]
\centering
\begin{tikzpicture}[x=1cm,y=.95cm,>=Stealth,
    every node/.style={font=\small},line join=round,
    edge/.style={draw=black,line width=.5pt}]
  \begin{scope}[shift={(.25,.25)},
      x={(2.15cm,0cm)},y={(.85cm,.55cm)},z={(0cm,1.9cm)}]
    \path[edge,fill=omniBlue!14] (0,0,0)--(1,0,0)--(1,0,.5)--(0,0,.5)--cycle;
    \path[edge,fill=omniOrange!18] (0,0,.5)--(.5,0,.5)--(.5,0,1)--(0,0,1)--cycle;
    \path[edge,fill=white] (.5,0,.5)--(1,0,.5)--(1,0,1)--(.5,0,1)--cycle;
    \path[edge,fill=omniBlue!14] (1,0,0)--(1,.5,0)--(1,.5,.5)--(1,0,.5)--cycle;
    \path[edge,fill=white] (1,0,.5)--(1,.5,.5)--(1,.5,1)--(1,0,1)--cycle;
    \path[edge,fill=black!22] (1,.5,0)--(1,1,0)--(1,1,1)--(1,.5,1)--cycle;
    \path[edge,fill=omniOrange!18] (0,0,1)--(.5,0,1)--(.5,1,1)--(0,1,1)--cycle;
    \path[edge,fill=white] (.5,0,1)--(1,0,1)--(1,.5,1)--(.5,.5,1)--cycle;
    \path[edge,fill=black!22] (.5,.5,1)--(1,.5,1)--(1,1,1)--(.5,1,1)--cycle;
    \node at (.5,0,.25) {$A$};
    \node at (.25,0,.75) {$B$};
    \node at (.75,0,.75) {$E$};
    \node at (1,.75,.5) {$C$};
    \draw[->] (0,0,-.15)--(.7,0,-.15) node[right] {$x$};
    \draw[->] (1.12,0,0)--(1.12,.75,0) node[right] {$y$};
    \draw[->] (-.12,0,0)--(-.12,0,.7) node[above] {$z$};
  \end{scope}
  \node[font=\footnotesize] at (1.75,3.05) {volume view};
  \begin{scope}[shift={(4.55,.45)}]
    \fill[omniBlue!14] (0,0) rectangle (1.7,.85);
    \fill[black!22] (.85,.85) rectangle (1.7,1.7);
    \draw[edge] (0,0) rectangle (1.7,1.7);
    \draw[edge] (0,.85)--(1.7,.85) (.85,.85)--(.85,1.7);
    \node at (.85,.425) {$A$};
    \node at (.425,1.275) {$D$};
    \node at (1.275,1.275) {$C$};
    \node[font=\footnotesize] at (.85,2) {$0\le z<\tfrac12$};
  \end{scope}
  \begin{scope}[shift={(7.15,.45)}]
    \fill[omniOrange!18] (0,0) rectangle (.85,1.7);
    \fill[black!22] (.85,.85) rectangle (1.7,1.7);
    \draw[edge] (0,0) rectangle (1.7,1.7);
    \draw[edge] (.85,0)--(.85,1.7) (.85,.85)--(1.7,.85);
    \node at (.425,.85) {$B$};
    \node at (1.275,.425) {$E$};
    \node at (1.275,1.275) {$C$};
    \node[font=\footnotesize] at (.85,2) {$\tfrac12\le z<1$};
  \end{scope}
  \draw[->] (4.3,.45)--(4.3,1.65) node[above] {$y$};
  \draw[->] (4.55,.2)--(5.75,.2) node[right] {$x$};
\end{tikzpicture}
\caption{The three-dimensional partition $\mathcal W$ and its two $z$-layers.
Boxes $A$, $B$, and $C$ span the $x$-, $y$-, and $z$-directions, respectively,
so every possible root cut crosses a leaf. All three views show
the same cell $C$. Unsplit extra coordinates give the same obstruction for $d>3$.}
\label{fig:inadmissible-partition}
\end{figure}

Some dyadic product partitions in $d\ge3$ are inadmissible
\cite{belkZaremskyTwisted2022}. For example, the disjoint boxes
\[
 [0,1)\!\times[0,\tfrac12)\!\times[0,\tfrac12),\quad
 [0,\tfrac12)\!\times[0,1)\!\times[\tfrac12,1),\quad
 [\tfrac12,1)\!\times[\tfrac12,1)\!\times[0,1)
\]
plus the two remaining half-size cubes partition $[0,1)^3$
(\Cref{fig:inadmissible-partition}). Each axis is spanned by one box,
so no root split is possible. We denote this partition
by $\mathcal W$ and its first box by $A$. Splitting $A$ in $x$ gives
an admissible partition with root mask $\{x\}$. Coarsening those siblings
back to $A$ would recover the inadmissible $\mathcal W$.
\begin{remark*}[Dimensions one and two]
In two dimensions every nontrivial dyadic partition has a common root split
axis. Otherwise a full-width and a full-height leaf would overlap,
hence coincide and equal the root cell. Recursion proves admissibility.
In one dimension this follows from interval nesting.
\end{remark*}

\subsubsection{Macro complexes}
\label{sec:macro-complexes}

A finite macro complex has reference hypercubes $[0,1]^d$ indexed by $\Mcal$,
each carrying one omnitree. Macro-face gluings couple these trees without physical
coordinates. Each face is paired reciprocally with a distinct face or carries
a boundary label \cite{mavriplisUnstructuredGridTechniques1997}.
Gluings permute and reverse tangential coordinates on $[0,1]^{d-1}$:
\begin{equation}
\label{eq:signed-permutation}
  \Phi(\zeta)_j=
  \begin{cases}
    \zeta_{\pi(j)},&\varepsilon_j=0,\\
    1-\zeta_{\pi(j)},&\varepsilon_j=1,
  \end{cases}
  \qquad j=1,\ldots,d-1,
\end{equation}
where \(\pi\) is a permutation of \([d-1]\) and \(\varepsilon_j\in\{0,1\}\). The
gluing acts on closed faces. \Cref{sec:trace-overlays} gives half-open representatives.

The forest $\Fcal=(\Tcal_m)_{m\in\Mcal}$ has active reference cells
$\Pcal(\Fcal):=\bigcup_{m\in\Mcal}\{m\}\times\Pcal(\Tcal_m)$.
A compatible physical embedding supplies the macro mesh of
\Cref{sec:finite-element-realization}.

\subsubsection{Cell identity and ordering}
\label{sec:cell-ordering}

A cell's identity is its key, the pair $(m,\kappa)$ of macro index and local
dyadic key. It is unchanged by tree regrouping or changes in stored position.
The active-key set is $\Lcal(\Fcal)=
\{(m,\kappa(v)):m\in\Mcal,\ v\in\leaves(\Tcal_m)\}$.
Its stored active-key array is sorted lexicographically by $(m,(\ell_a)_a,(i_a)_a)$.
For a cell $K$, $\ell_a(K),i_a(K)$ are its key components and
$\ch_\sigma(K)$ its child cells.
Lexicographic order groups macro leaves. Preorder groups subtrees.

Every internal node has at least two children, so a tree with $N_m$ leaves
has at most $2N_m-1$ records.

\subsection{Computational model}
\label{sec:computational-model}
We count sequential word-RAM work on flat storage: contiguous arrays of
forest records, requests and materialized outputs. \emph{Output-optimal} work is proportional
to the explicitly stored input and selected output. A \emph{handle} identifies
an immutable stored subtree. Reduction can reuse such handles internally.
The final forest and selected hierarchy levels are emitted as complete arrays.
For a fixed key width $w\ge L_{\max}$, operations on $w$-bit words have unit cost.
Key ordering uses $O_d(w)$ radix passes.
Macro identifiers, descriptor handles and all array indices fit one word,
including intermediate and output indices. Keys use $O(d)$ words.
The bounds retain level dependence: $O_d$ hides only dimension, whereas
$O_{d,w}$ also fixes key width. Supplied offsets give direct access. Other lookup, ordering,
and branch encoding costs are charged.

\Needspace{8\baselineskip}
\section{Canonical representation}
\label{sec:active-normal}
\label{sec:normal-form}

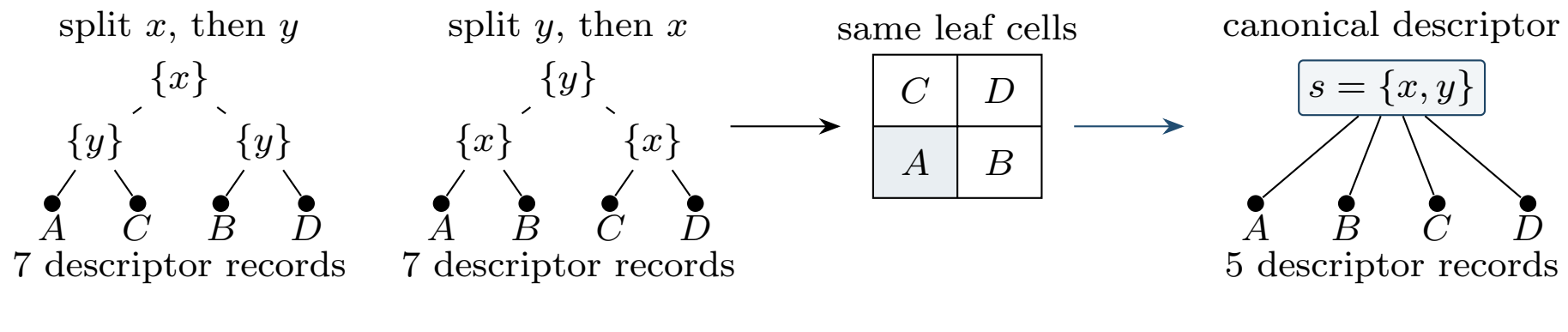
\begin{figure}[!hbp]
\centering
\begin{tikzpicture}[x=1cm,y=.85cm,>=Stealth,
    every node/.style={font=\footnotesize},
    mask/.style={inner sep=2pt,fill=white},
    tip/.style={circle,fill=black,inner sep=1.3pt}]
  \foreach \xx/\first/\second/\names in {0/x/y/{A,C,B,D},3/y/x/{A,B,C,D}} {
    \begin{scope}[shift={(\xx,0)}]
      \node at (0,1.6) {split $\first$, then $\second$};
      \node[mask] (r\first) at (0,1.13) {$\{\first\}$};
      \foreach \u/\v in {0/-.65,1/.65} {
        \node[mask] (n\first\u) at (\v,.56) {$\{\second\}$};
        \draw (r\first)--(n\first\u);
      }
      \foreach \k/\v in {0/-.98,1/-.32,2/.32,3/.98} {
        \node[tip] (t\first\k) at (\v,0) {};
        \pgfmathtruncatemacro{\parent}{floor(\k/2)}
        \draw (n\first\parent)--(t\first\k);
      }
      \foreach \name [count=\k from 0] in \names {
        \pgfmathsetmacro{\v}{-.98+.653333*\k}
        \node at (\v,-.22) {$\name$};
      }
      \node at (0,-.58) {7 descriptor records};
    \end{scope}
  }
  \node at (6,1.6) {same leaf cells};
  \begin{scope}[shift={(5.35,.05)}]
    \fill[omniBlue!10] (0,0) rectangle (.65,.65);
    \draw (0,0) rectangle (1.3,1.3);
    \draw (.65,0)--(.65,1.3) (0,.65)--(1.3,.65);
    \foreach \x/\y/\name in {.325/.325/A,.975/.325/B,.325/.975/C,.975/.975/D}
      \node at (\x,\y) {$\name$};
  \end{scope}
  \draw[->] (4.25,.7)--(5.1,.7);
  \draw[->,omniBlue!90!black] (6.9,.7)--(7.75,.7);
  \node at (9.35,1.6) {canonical descriptor};
  \node[mask,draw=omniBlue!80!black,rounded corners=1pt,fill=omniBlue!6]
    (canon) at (9.35,1.05) {$s=\{x,y\}$};
  \foreach \v/\name in {8.3/A,9/B,9.7/C,10.4/D} {
    \node[tip] (leaf\name) at (\v,0) {};
    \draw (canon)--(leaf\name);
    \node at (\v,-.22) {$\name$};
  }
  \node at (9.35,-.58) {5 descriptor records};
\end{tikzpicture}
\caption{Lifting a common split axis (\Cref{lem:axis-lifting}) gives one
canonical descriptor for both histories. Dots are leaf records.}
\label{fig:normal-form-commutation}
\end{figure}

The four quadrants of a square can be produced by splitting $x$ then $y$,
$y$ then $x$, or both axes together (\Cref{fig:normal-form-commutation}).
Only the tree grouping differs. Canonical reconstruction chooses one grouping:
at each nonleaf region, it takes every midpoint cut that crosses no leaf
interior, then recurses on the children. It preserves the cells and succeeds
exactly when their partition has a tree representation
(\Cref{thm:normal-form}).

$\Lcal(U;\kappa)$ is the set of leaf keys of $U$ rooted at $\kappa$;
we omit $\kappa$ when the root is known. The admissible leaf-key sets are
$\Acal(\kappa)=\{\Lcal(U;\kappa):U\text{ valid at }\kappa\}$.
For leaf keys $\Lcal$ partitioning $Q_\kappa$, the \emph{maximal common mask} is
\begin{equation}
\label{eq:root-mask}
  s(\kappa,\Lcal)=
  \set{a\in[d]\given
    \ell_a(\eta)>\ell_a(\kappa)\text{ for every }\eta\in\Lcal}.
\end{equation}
For $a\in s(\kappa,\Lcal)$, each leaf lies in one of
$Q_{c_{\{a\},0}(\kappa)}$ and $Q_{c_{\{a\},1}(\kappa)}$ by
\eqref{eq:containment}. An axis outside $s$ has a leaf spanning the region
in that axis, preventing a root cut. For
\(s=s(\kappa,\Lcal)\ne\emptyset\) and \(b\in\{0,1\}^s\), put
\begin{equation}
\label{eq:child-leaf-set-formal}
  \Lcal_b:=\set{\eta\in\Lcal\given c_{s,b}(\kappa)\preceq\eta} .
\end{equation}

If the root mask of \(U\) contains \(a\), write \(U[a=\varepsilon]\) for the
descriptor obtained by selecting the branches with \(a\)-bit \(\varepsilon\)
and removing \(a\) from the root mask. If the remaining mask is empty,
suppress the root. This retains the leaf keys contained in
\(Q_{c_{\{a\},\varepsilon}(\kappa)}\).

The binary form of \Cref{lem:axis-lifting} is proved in
\cite{hennigMatucciPresentations2012}*{Lemma~6}. We give the simultaneous-mask
proof used by the admissibility criterion in \Cref{thm:normal-form}.

\begin{lemma}[A common split axis can be moved to the root]
\label{lem:axis-lifting}
Let $\Lcal\in\Acal(\kappa)$, $\Lcal\ne\{\kappa\}$, and
$a\in s(\kappa,\Lcal)$. Given a descriptor for $\Lcal$ with root mask $\sigma$, there is a
descriptor with root mask $\sigma\cup\{a\}$.
\end{lemma}

\begin{proof}
Induct on $N=|\Lcal|$. For $N=1$ there is no axis to lift; $a\in\sigma$
needs no change. Otherwise set $\kappa_b=c_{\sigma,b}(\kappa)$ and let
$\Lcal_b^\sigma$ be its leaf-key set. Then
$\Lcal_b^\sigma\in\Acal(\kappa_b)$, $|\Lcal_b^\sigma|<N$, and
$\ell_a(\kappa_b)=\ell_a(\kappa)
<\min_{\eta\in\Lcal_b^\sigma}\ell_a(\eta)$.
Thus $a\in s(\kappa_b,\Lcal_b^\sigma)$, and induction gives descriptors
$U_b$ whose roots contain $a$. Set
\[
 \widehat U=\bigl(\sigma\cup\{a\},(U_b[a=\varepsilon])_{b,\varepsilon}\bigr),
 \qquad
 c_{\{a\},\varepsilon}(\kappa_b)
 =c_{\sigma\cup\{a\},(b,\varepsilon)}(\kappa).
\]
By \eqref{eq:child}, the child keys agree, and
$\Lcal(\widehat U;\kappa)
=\dotcup_{b,\varepsilon}\Lcal(U_b[a=\varepsilon])=\Lcal$.
\end{proof}

\subsubsection*{Canonical reconstruction}
For a nontrivial admissible partition, successive lifts give root mask $s$,
since every descriptor for $\Lcal$ has
$\emptyset\ne\sigma\subseteq s$. Its children satisfy
\begin{equation}
\label{eq:leaf-disjoint-recursion}
  \Lcal=\dotcup_{b\in\{0,1\}^{s}}\Lcal_b,
  \qquad
  \emptyset\ne \Lcal_b\in\Acal(c_{s,b}(\kappa)).
\end{equation}

Canonical reconstruction is defined recursively by
\begin{equation}
\label{eq:normal-recursion}
\Reconstruct_{\kappa}(\Lcal)=
\begin{cases}
  \leaf,&\Lcal=\{\kappa\},\\[1mm]
  \bigl(s,(\Reconstruct_{c_{s,b}(\kappa)}(\Lcal_b))_{b\in\{0,1\}^s}\bigr),
  &\text{otherwise},
\end{cases}
\end{equation}
where $s=s(\kappa,\Lcal)$, in the branch order of \Cref{def:descriptor}.

On any dyadic partition, the recursion \eqref{eq:normal-recursion} fails if a
nonsingleton subpartition has $s=\emptyset$. Otherwise nesting assigns each leaf
to one child and coverage makes all child sets nonempty.
Canonicalization is $\operatorname{Can}_{\kappa}(U)
=\Reconstruct_{\kappa}(\Lcal(U;\kappa))$. Below, $U_v$ denotes a node's subtree.

\begin{theorem}[Admissibility and canonical representation]
\label{thm:normal-form}
A finite dyadic partition of $Q_{\kappa}$ is admissible if and only if
canonical reconstruction \eqref{eq:normal-recursion} succeeds. For $\Lcal\in\Acal(\kappa)$,
$\Reconstruct_{\kappa}(\Lcal)$ is the unique valid descriptor with leaf-key set
$\Lcal$ satisfying
\begin{equation}
\label{eq:maximal-mask-characterization}
  \operatorname{rm}(U_v)=s(\kappa(v),\Lcal(U_v;\kappa(v)))
  \quad\text{at every internal node }v.
\end{equation}
Thus $\operatorname{Can}_{\kappa}(U)=\operatorname{Can}_{\kappa}(U')$
if and only if $\Lcal(U;\kappa)=\Lcal(U';\kappa)$.
\end{theorem}

\begin{proof}
Induct on $N=|\Lcal|$. For $N=1$, coverage gives $\Lcal=\{\kappa\}$,
represented only by $\leaf$.
For $N>1$, a valid descriptor has a nonempty root mask
$\sigma\subseteq s$. By \Cref{lem:axis-lifting} and \eqref{eq:leaf-disjoint-recursion},
each $\Lcal_b$ is admissible with $1\le|\Lcal_b|<N$.
Induction constructs every child in \eqref{eq:normal-recursion}.
Conversely, successful recursive calls and \eqref{eq:child-partition} assemble a
valid descriptor. Every call has fewer leaves, so the attempt terminates.

For uniqueness, \eqref{eq:maximal-mask-characterization} forces root mask
$s$ and child leaf-key sets $\Lcal_b$. Induction forces each $U_b$, hence $U$.
Equal leaf-key sets give equal canonical descriptors. Decoding gives the converse.
\end{proof}

A descriptor is \emph{canonical} if $\operatorname{Can}_\kappa(U)=U$.
\Cref{thm:normal-form} gives $\operatorname{Can}^2=\operatorname{Can}$.
By induction, \eqref{eq:maximal-mask-characterization} is equivalent to the
local normalization condition stipulated for uniqueness in
\cite{pollingerBeautyAnisotropic2026}*{Definition~2.2}. Here we prove uniqueness
together with the admissibility criterion and give the flat construction below.
\begin{proposition}[Canonical reconstruction from a leaf-key array]
\label{prop:flat-normalization}
Given an array $(\eta_1,\ldots,\eta_N)$ of admissible leaf keys with coordinate
levels at most $L_{\max}$, direct evaluation of \eqref{eq:normal-recursion} constructs
its canonical descriptor in $O_d((1+L_{\max})N)$ work and $O_d(N)$ workspace,
with at most $2N-1$ descriptor records. Thus \Cref{sec:computational-model}
gives $O_{d,w}(N)$ work.
\end{proposition}
\begin{proof}
At each visited region $\kappa$, scan its leaf range $\Lcal_\kappa$,
bucket by child, and recurse. For the visited regions $\mathcal R$,
\begin{equation}
\label{eq:flat-reconstruction-charge}
\begin{aligned}
 \sum_{\kappa\in\mathcal R}|\Lcal_\kappa|
 &=\sum_{\eta\in\Lcal}\card{\{\kappa\in\mathcal R:\kappa\preceq\eta\}}\\
 &\le\sum_{\eta\in\Lcal}\left(1+\sum_a\ell_a(\eta)\right)
 \le(1+dL_{\max})N.
\end{aligned}
\end{equation}
This bounds scans and bucketing. Two flat key buffers give $O_d(N)$
workspace; \Cref{sec:cell-ordering} gives at most $2N-1$ records.
\end{proof}

\section{Admissible directional adaptation}
\label{sec:normalization-adaptation}

Our adaptation preserves a tree representation. A \emph{contraction} replaces
two permitted dyadic siblings by their rectangular union, undoing one coordinate
split. Canonical regrouping changes only the tree and can
expose further permitted pairs (\Cref{fig:directional-coarsening}).
Each tree is adapted independently. Face closure then couples neighboring
cells across the forest.

A \emph{valid one-level plan} assigns refinement masks $R(\eta)$ and
contraction permission masks $C(\eta)$ to input active keys. Submitted masks are
nonempty. Duplicates are united and unrequested masks are $\emptyset$.
Each input leaf $\eta$ is refined once in $R(\eta)$, if nonempty. Its children
\eqref{eq:child} must lie in $\Kcal_d$ \eqref{eq:dyadic-key-cell}.
Refinement takes precedence: children receive $\chi=C(\eta)\setminus R(\eta)$. Unrefined leaves
retain $C(\eta)$. Permission masks allow, but do not require, contraction.
The one-level restriction applies to adaptation. Subsequent face closure may
refine further.

\subsection{Permitted contractions}

An exposed pair consists of two leaf children of the same node that differ
in one branch bit. Contracting the pair requires both leaves to permit that axis.
At a root $(\sigma,(V_b)_b)$, write $\leaf_\chi$ for a leaf with permission mask $\chi$.
For $a\in\sigma$, write $r_\varepsilon=r\cup\{a\mapsto\varepsilon\}$.
An \emph{exposed permitted $a$-pair} means
\begin{equation}
\label{eq:exposed-pair}
 r\in\{0,1\}^{\sigma\setminus\{a\}},\qquad
 V_{r_0}=\leaf_{\chi_0},\quad V_{r_1}=\leaf_{\chi_1},\qquad
 a\in\chi_0\cap\chi_1.
\end{equation}
Geometric siblings with different parents are \emph{hidden}: contracting them could recover the
inadmissible $\mathcal W$ (\Cref{fig:inadmissible-partition}).
Such pairs can be contracted only after canonical regrouping exposes them.
A contraction replaces the pair by $\leaf_{(\chi_0\cap\chi_1)\setminus\{a\}}$.
\subsection{Canonical regrouping and directional reduction}
\label{sec:descriptor-normalizer}
\label{lem:section-fusion}

Section selects root branches without cutting leaves. Fusion joins adjacent
subtrees under a canonical root, lifting additional common split axes to the
parent (\Cref{fig:section-fusion,fig:normal-form-commutation}). Both preserve
leaf cells and permission masks.
For branch bits, $c|_\sigma$ restricts $c$ to $\sigma$, while $e\cup b$ combines
assignments on disjoint axes.

For a canonical descriptor $V$ rooted at $\kappa$ with root mask $\sigma$,
$f\subseteq\sigma$, and $e\in\{0,1\}^f$, let
$\operatorname{Sec}(V;f=e)$ canonically represent the leaves selected by
the original root bits $f=e$, rooted at $c_{f,e}(\kappa)$ when $f\ne\emptyset$.
Then
\begin{equation}
\label{eq:canonical-section}
  \operatorname{Sec}(V;f=e)=
  \begin{cases}
    V,&f=\emptyset,\\
    V_e,&\emptyset\ne f=\sigma,\\
    \operatorname{Fuse}_{\sigma\setminus f}
      \bigl((V_{e\cup b})_{b\in\{0,1\}^{\sigma\setminus f}}\bigr),
      &\emptyset\ne f\subsetneq\sigma.
  \end{cases}
\end{equation}

For canonical child descriptors $V_b$ on $c_{\sigma,b}(\kappa)$ with leaf-key sets $\Lcal_b$,
$\sigma\ne\emptyset$ lists the parent's existing splits. Fusion lifts additional
axes $\alpha$ common to all children, giving canonical root mask $s=\sigma\cup\alpha$:
\begin{equation}
\label{eq:canonical-fusion}
\begin{gathered}
  \alpha=\left(\bigcap_{b\in\{0,1\}^{\sigma}}\operatorname{rm}(V_b)\right)\setminus\sigma,
  \qquad s=\sigma\cup\alpha,\\
  \operatorname{Fuse}_{\sigma}((V_b)_b):=
  \left(s,
    \left(\operatorname{Sec}
      (V_{c|_{\sigma}};\alpha=c|_\alpha)\right)_{c\in\{0,1\}^s}\right).
\end{gathered}
\end{equation}
Each nonidentity recursive call has fewer participating leaves, proving
termination. By \eqref{eq:root-mask}, canonical child descriptors and equal parent/child levels outside $\sigma$ imply
$\sigma\cup\alpha=s(\kappa,\bigcup_b\Lcal_b)$.
Induction on participating leaves, using
\eqref{eq:maximal-mask-characterization}, proves canonicality.
\Cref{lem:section-work} bounds work.

\subsubsection*{Directional reduction}
Permission masks limit coarsening along a \emph{point lineage}: the successive
leaves containing a fixed reference point.

\Needspace{8\baselineskip}
\begin{lemma}[One contraction per axis along a point lineage]
\label{lem:lineage-invariant}
During one adaptation, each axis is contracted at most once along every
point lineage.
Along the same lineage, an axis refined in the initial step cannot
subsequently be contracted.
At a contraction, the resulting permission mask excludes every axis
contracted in either child's history.
\end{lemma}
\begin{proof}
Let $H_i$ contain axes initially refined or already contracted in a leaf's
history. Initially $\chi_i\cap H_i=\emptyset$. At an $a$-contraction,
\begin{equation}
\label{eq:permission-history}
 H'=H_0\cup H_1\cup\{a\},\qquad
 \chi'=(\chi_0\cap\chi_1)\setminus\{a\},\qquad
 \chi'\cap H'=\emptyset.
\end{equation}
Section and fusion preserve $(\chi,H)$, while
$a\in\chi_0\cap\chi_1$ implies $a\notin H_0\cup H_1$.
Induction excludes both repeated and initially refined axes.
\end{proof}

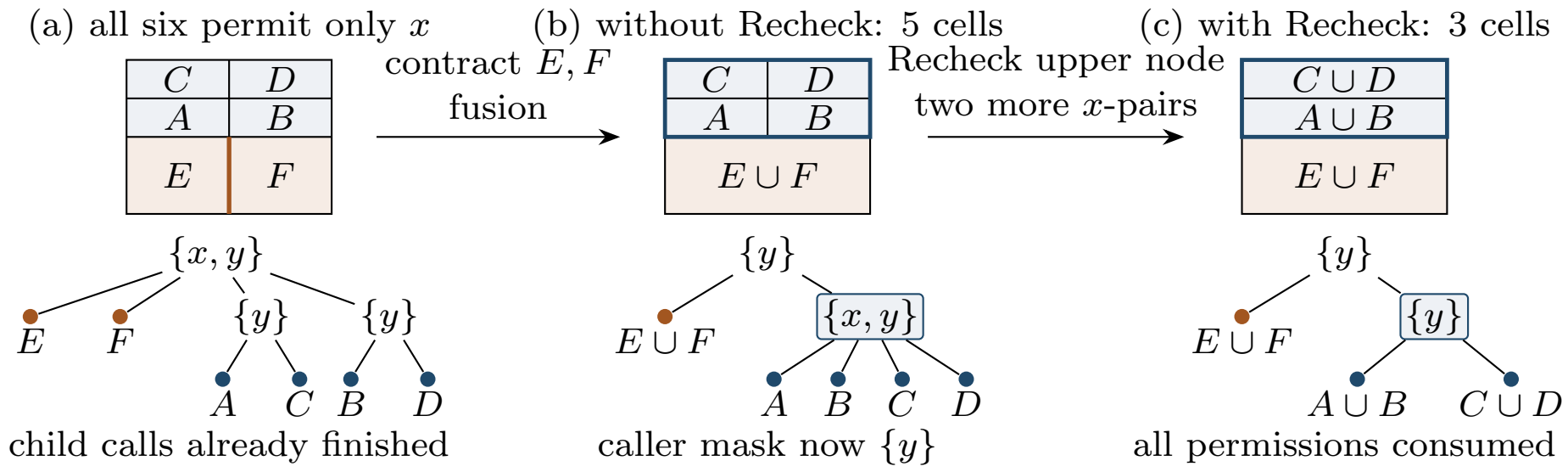
\begin{figure}[!htbp]
\centering
\begin{tikzpicture}[x=1cm,y=.75cm,>=Stealth,
  every node/.style={font=\footnotesize},
  mask/.style={fill=white,inner sep=1pt},
  tip/.style={circle,inner sep=1.2pt},
  newroot/.style={mask,draw=omniBlue!85!black,rounded corners=1pt,fill=omniBlue!8}]
  \foreach \xx/\stage/\title in {.8/0/{(a) all six permit only $x$},5/1/{(b) without Recheck: 5 cells},9.5/2/{(c) with Recheck: 3 cells}} {
    \begin{scope}[shift={(\xx,0)}]
      \node at (.8,2.25) {\title};
      \fill[omniOrange!12] (0,.3) rectangle (1.6,1.1);
      \fill[omniBlue!8] (0,1.1) rectangle (1.6,1.9);
      \draw (0,.3) rectangle (1.6,1.9);
      \draw (0,1.1)--(1.6,1.1) (0,1.5)--(1.6,1.5);
      \ifnum\stage<2
        \draw (.8,1.1)--(.8,1.9);
        \foreach \x/\low/\high in {.4/A/C,1.2/B/D} {
          \node at (\x,1.3) {$\low$};
          \node at (\x,1.7) {$\high$};
        }
      \else
        \node at (.8,1.3) {$A\cup B$};
        \node at (.8,1.7) {$C\cup D$};
      \fi
      \ifnum\stage=0
        \draw[omniOrange!90!black,line width=1pt] (.8,.3)--(.8,1.1);
        \node at (.4,.7) {$E$}; \node at (1.2,.7) {$F$};
      \else
        \node at (.8,.7) {$E\cup F$};
        \draw[omniBlue!85!black,line width=.9pt] (0,1.1) rectangle (1.6,1.9);
      \fi
    \end{scope}
  }
  \draw[->] (2.75,1.1)--(4.65,1.1);
  \node[align=center] at (3.7,1.65) {contract $E,F$\\fusion};
  \draw[->] (7.05,1.1)--(9.05,1.1);
  \node[align=center] at (8.05,1.65) {Recheck upper node\\two more $x$-pairs};
  \begin{scope}[yshift=2.5mm]
  \node[mask] (r0) at (1.5,-.45) {$\{x,y\}$};
  \node[tip,fill=omniOrange!90!black] (E0) at (.05,-1.1) {};
  \node[tip,fill=omniOrange!90!black] (F0) at (.75,-1.1) {};
  \node[mask] (l0) at (1.85,-1.1) {$\{y\}$};
  \node[mask] (u0) at (2.85,-1.1) {$\{y\}$};
  \foreach \p in {E0,F0,l0,u0} \draw (r0)--(\p);
  \node at (.05,-1.35) {$E$}; \node at (.75,-1.35) {$F$};
  \foreach \x/\name/\p in {1.55/A/l0,2.15/C/l0,2.55/B/u0,3.15/D/u0} {
    \node[tip,fill=omniBlue!85!black] (i\name) at (\x,-1.75) {};
    \draw (\p)--(i\name);
    \node at (\x,-2) {$\name$};
  }
  \node at (1.6,-2.45) {child calls already finished};
  \node[mask] (r1) at (5.8,-.45) {$\{y\}$};
  \node[tip,fill=omniOrange!90!black] (s1) at (5,-1.1) {};
  \node[newroot] (t1) at (6.6,-1.1) {$\{x,y\}$};
  \draw (r1)--(s1) (r1)--(t1);
  \node at (5,-1.35) {$E\cup F$};
  \foreach \x/\name in {5.85/A,6.35/B,6.85/C,7.35/D} {
    \node[tip,fill=omniBlue!85!black] (m\name) at (\x,-1.75) {};
    \draw (t1)--(m\name);
    \node at (\x,-2) {$\name$};
  }
  \node at (5.8,-2.45) {caller mask now $\{y\}$};
  \node[mask] (r2) at (10.3,-.45) {$\{y\}$};
  \node[tip,fill=omniOrange!90!black] (s2) at (9.5,-1.1) {};
  \node[newroot] (t2) at (11,-1.1) {$\{y\}$};
  \draw (r2)--(s2) (r2)--(t2);
  \node at (9.5,-1.35) {$E\cup F$};
  \foreach \x/\name/\cell in {10.4/ab/{A\cup B},11.6/cd/{C\cup D}} {
    \node[tip,fill=omniBlue!85!black] (o\name) at (\x,-1.75) {};
    \draw (t2)--(o\name);
    \node at (\x,-2) {$\cell$};
  }
  \node at (10.3,-2.45) {all permissions consumed};
  \end{scope}
\end{tikzpicture}
\caption{All six leaves permit only $x$-contraction. After $E,F$ contract,
fusion creates the blue root. Recheck calls Merge there
(line~\ref{line:recheck-merge}), giving three cells instead of five.}
\label{fig:directional-coarsening}
\label{fig:section-fusion}
\end{figure}

Let $\Lcal_R$ be the post-refinement leaf-key set and
$V_R=\Reconstruct_{\rootkey}(\Lcal_R)$ its canonical descriptor
\eqref{eq:normal-recursion}, with contraction permissions $\chi$ from
\Cref{sec:normalization-adaptation}.
Reduce processes children first. Merge chooses the root axis with the most
exposed permitted pairs, contracts them and fuses the other slots. Recheck
tests their new roots without repeating child reductions.
With only $x$-permissions in \Cref{fig:directional-coarsening}, $c_x=1$, $c_y=0$, so
$E,F$ contract first. Fusion creates an upper $\{x,y\}$ root while the
caller retains $\{y\}$. Recheck contracts $A,B$ and $C,D$ there.
If only $A,B$ have permission $\{x\}$, no exposed pair can contract and all six
cells remain. Thus the policy can leave an admissible permitted merge undone.

\begin{algorithm}[htbp]
\caption{Directional reduction}
\label{alg:directional-reduction}
{\small
\begin{algorithmic}[1]
\Function{Reduce}{$U$}
  \State \textbf{if} $U$ is a leaf \textbf{then return} $U$
  \State $(V_b)_b\gets(\Call{Reduce}{U_b})_b$ in branch order \label{line:reduce-children}
  \State \Return $\Call{Merge}{\operatorname{rm}(U),(V_b)_b}$ \label{line:reduce-merge}
\EndFunction
\Function{Merge}{$\sigma,(V_b)_b$}
  \State \textbf{if} $\sigma=\emptyset$ \textbf{then return} the sole member
  \State $c_a\gets$ number of exposed permitted $a$-pairs, for each $a\in\sigma$
  \State \textbf{if} $\max_{a\in\sigma}c_a=0$ \textbf{then return}
    $\operatorname{Fuse}_\sigma((V_b)_b)$ \label{line:merge-fuse}
  \State choose $a$ maximizing $c_a$, ties to the smallest axis \label{line:merge-choice}
  \For{each slot $r\in\{0,1\}^{\sigma\setminus\{a\}}$ in branch order}
    \If{$V_{r_0},V_{r_1}$ are leaves and both permit $a$}
      \State $W_r\gets$ a leaf on their union with permission mask
      $ (\chi_0\cap\chi_1)\setminus\{a\}$ \label{line:merge-contract}
    \Else
      \State $W_r\gets\Call{Recheck}{\operatorname{Fuse}_{\{a\}}
        ((V_{r_0},V_{r_1}))}$ \label{line:merge-recheck}
    \EndIf
  \EndFor
  \State \Return $\Call{Merge}{\sigma\setminus\{a\},(W_r)_r}$ \label{line:merge-residual}
\EndFunction
\Function{Recheck}{$V$}
  \State \textbf{if} $V$ is a leaf \textbf{then return} $V$
  \State \Return $\Call{Merge}{\operatorname{rm}(V),
    \text{canonical root children of }V}$ \label{line:recheck-merge}
\EndFunction
\end{algorithmic}
}
\end{algorithm}

Write $\Apply(U;R,C)=\Lcal(\operatorname{Reduce}(V_R))$,
bypassing Reduce when all $C(\eta)$ are empty. Refinement-only adaptation may pass
flat leaves to closure and defer canonical reconstruction to final emission.
For initial face matching, a working copy of the input tree with requested
splits grafted onto its leaves and child locators suffices.
This takes $O_d(n)$ work for $n$ input descriptor records. Matching requires no
canonical regrouping (\Cref{sec:interface-output-complexity}).
\begin{proposition}[Admissible directional reduction]
\label{thm:adaptation-correctness}
\label{lem:admissible-contraction}
Each exposed contraction \eqref{eq:exposed-pair} preserves admissibility.
For a valid input descriptor $U$ and a valid one-level plan,
Algorithm~\ref{alg:directional-reduction}, applied to the canonical
post-refinement descriptor $V_R$, terminates and returns deterministic
$U'$ satisfying
\[
 \Apply(U;R,C)\in\Acal(\rootkey),\qquad
 U'=\Reconstruct_{\rootkey}(\Apply(U;R,C)).
\]
\end{proposition}
\begin{proof}
\emph{Admissibility.} Grafting the requested splits onto $U$ preserves
admissibility, so $V_R$ exists by \Cref{thm:normal-form}.
For an exposed contraction, put $\sigma'=\sigma\setminus\{a\}$. For each residual slot $r$,
let $W_r$ be a leaf on the union of a contracted pair, or retain
$W_r=(\{a\},(V_{r_0},V_{r_1}))$ otherwise. Then
$\widetilde U=(\sigma',(W_r)_r)$, with an empty root suppressed, certifies
admissibility by \eqref{eq:child-partition}. Section and fusion preserve
these leaves by \eqref{eq:canonical-section}--\eqref{eq:canonical-fusion}.

\emph{Termination.} Only lines~\ref{line:reduce-merge} and \ref{line:recheck-merge}
can have $k'=k$, with $F'=\operatorname{Merge}\ne F$.
All other calls have $k'<k$: child regions are proper, a residual call
follows a contraction, and a failed slot has fewer leaves since a successful
slot exists.
For $M(F,k)=2k+\1_{\{F\ne\operatorname{Merge}\}}$,
\begin{equation}
\label{eq:reduction-measure}
 M(F',k')\le
 \begin{cases}
  2k<2k+1=M(F,k),&k'=k,\\
  2k'+1<2k\le M(F,k),&k'<k.
 \end{cases}
\end{equation}
Section and fusion terminate by \Cref{sec:descriptor-normalizer}.

\emph{Canonicality.} Induction on $M$ using \eqref{eq:canonical-fusion}
proves canonical returns. Branch order and line~\ref{line:merge-choice}
determine the partition. \Cref{thm:normal-form} gives uniqueness.
\end{proof}

\subsection{Adaptation cost}
\label{sec:adaptation-cost}

Regrouping can place the same dyadic region under different tree parents.
Two counts control its work: how often a region is visited, and how many
regions can contain one post-refinement leaf. The lemma bounds the first.
The proposition then sums the second.
Let $n,n_{\rm ad}$ count input and adapted descriptor records. One-level
refinement creates at most $2^dn$ leaves, whose canonical descriptor is
constructed by \Cref{prop:flat-normalization}. Contraction decreases this
count, so $n_{\rm ad}=O_d(n)$.

\Needspace{10\baselineskip}
\begin{lemma}[Region-event budget]
\label{lem:region-events}
Fix a region $\kappa$ in one execution of Algorithm~\ref{alg:directional-reduction}
on $V_R$.
A Merge call is \emph{successful} if it directly contracts an exposed pair.
Let $s_\kappa$ count successes, $e_\kappa$ Reduce/Recheck entries, and
$m_\kappa$ all Merge calls. For the one-step product parents
$\Pi_\kappa=\{\nu:\kappa=c_{\sigma,b}(\nu),\ \sigma\ne\emptyset\}$,
\begin{equation}
\label{eq:reduction-call-budget}
 s_\kappa\le d,\qquad e_\kappa\le1+\sum_{\nu\in\Pi_\kappa}s_\nu,
 \qquad m_\kappa\le e_\kappa+s_\kappa\le1+d2^d.
\end{equation}
\end{lemma}
\begin{proof}
A successful $a$-call creates a leaf spanning $\kappa$'s $a$-interval.
Later coarsening contains this witness. Section and fusion preserve leaves.
Every later call at $\kappa$ therefore excludes $a$ from its split mask.
A witness enlarged outside $\kappa$ prevents further calls there.
Thus $s_\kappa\le d$, attained by a full tensor split with full permissions.

Reduce enters an original region at most once. Each Recheck entry
(line~\ref{line:merge-recheck}) comes from a successful one-step parent call,
with at most one entry per parent success; $|\sigma|=1$ creates no Recheck.
Each residual Merge (line~\ref{line:merge-residual}) follows a success.
Each nonempty split mask gives at most one product parent. Consequently,
$|\Pi_\kappa|\le2^d-1$ and
$m_\kappa\le e_\kappa+s_\kappa\le1+d(2^d-1)+d=1+d2^d$.
\end{proof}

\begin{proposition}[Work of directional reduction]
\label{prop:reduction-work}
Let $\Lcal_R=\Lcal(V_R)$ contain $N$ post-refinement leaves. With immutable
descriptor handles and child locators,
\[
 T_{\rm reduce}+S_{\rm reduce}\le C_d\sum_{\eta\in\Lcal_R}
 \left(1+\frac12\sum_{a=1}^d\ell_a(\eta)\right)
 \prod_{a=1}^d(\ell_a(\eta)+1),
\]
where $S_{\rm reduce}$ counts allocated descriptor records and $C_d$ depends
only on dimension. This is $O_d((1+L)^{d+1}N)$ for maximal coordinate level $L$.
After request resolution, one-level adaptation therefore takes
$O_{d,w}(n)$ work and storage, including directional coarsening.
\end{proposition}
\begin{proof}
Let $\mathcal R$ be the distinct call roots and set
$W_\kappa=\sum_{\substack{\eta\in\Lcal_R\\\kappa\preceq\eta}}q_\kappa(\eta)$, where
$q_\kappa(\eta)=1+\|\ell(\eta)-\ell(\kappa)\|_1$.

For the current leaves $\zeta$ at root $\kappa$, their groups of contained
post-refinement leaves are disjoint and nonempty. Since
$\kappa\preceq\zeta\preceq\eta$ implies
$q_\kappa(\zeta)\le q_\kappa(\eta)$, each call's local work and allocation
are bounded by $C_dW_\kappa$ (\Cref{lem:section-work}).
Call multiplicities are bounded by \Cref{lem:region-events}. Each product
ancestor is fixed by its levels $k_a$ \eqref{eq:containment}. Thus
\begin{equation}
\label{eq:reduction-ancestor-charge}
\begin{aligned}
 T_{\rm reduce}+S_{\rm reduce}
 &\le C_d\sum_{\kappa\in\mathcal R}(e_\kappa+m_\kappa)W_\kappa
 \le C_d'\sum_{\kappa\in\mathcal R}W_\kappa\\
 &\le C_d'\sum_{\eta\in\Lcal_R}
 \sum_{0\le k_a\le\ell_a(\eta)}
 \left(1+\sum_a(\ell_a(\eta)-k_a)\right)\\
 &=C_d'\sum_{\eta\in\Lcal_R}
 \left(1+\frac12\sum_a\ell_a(\eta)\right)
 \prod_a(\ell_a(\eta)+1).
\end{aligned}
\end{equation}
If only $r$ coordinates grow in depth, the factor is
$O((1+L)^{r+1})$. Reconstruction costs $O_d((1+L)N)$
(\Cref{prop:flat-normalization}), with $N=O_d(n)$ and $L\le w$.
\end{proof}

\section{Multilevel hierarchies and cell incidence}
\label{sec:multilevel}
\label{sec:hierarchy}

\subsection{Rank cuts and cell incidence}
\label{sec:transfer}

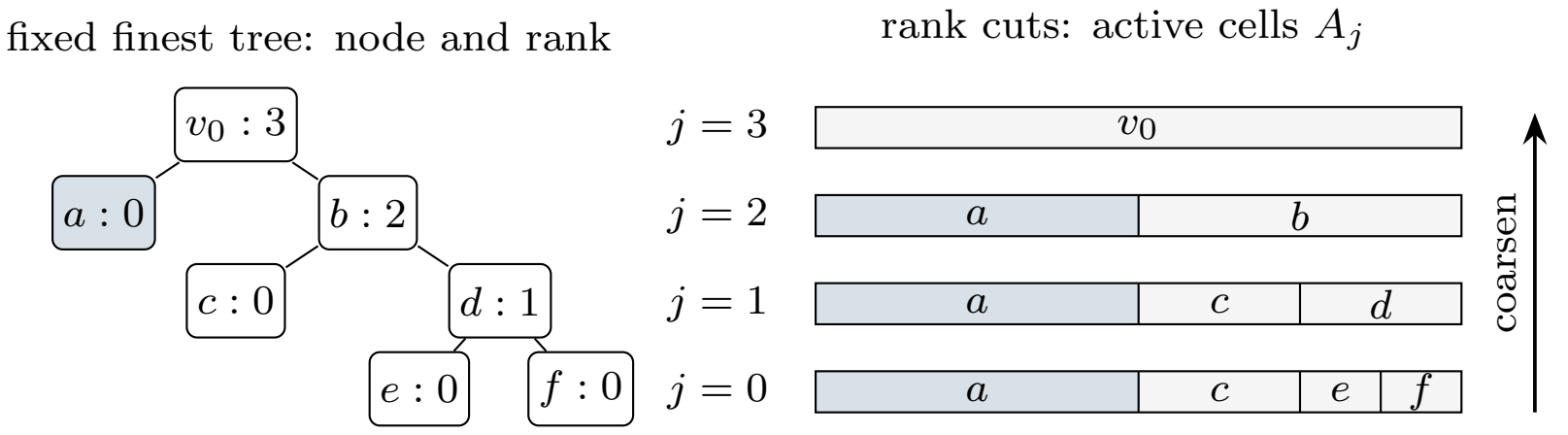
\begin{figure}[!htbp]
\centering
\begin{tikzpicture}[x=1cm,y=.8cm,>=Stealth,
    every node/.style={font=\footnotesize},
    rank/.style={draw,rounded corners=2pt,fill=white,inner sep=2pt,minimum height=5mm}]
  \begin{scope}[shift={(-.25,0)}]
  \node[anchor=south] at (1.9,.83) {fixed finest tree: node and rank};
  \node[rank] (root) at (1.4,.35) {$v_0:3$};
  \node[rank,fill=omniBlue!18] (left) at (.5,-.4) {$a:0$};
  \node[rank] (right) at (2.3,-.4) {$b:2$};
  \node[rank] (middle) at (1.4,-1.15) {$c:0$};
  \node[rank] (deep) at (3.2,-1.15) {$d:1$};
  \node[rank] (lastleft) at (2.65,-1.9) {$e:0$};
  \node[rank] (lastright) at (3.75,-1.9) {$f:0$};
  \foreach \parent/\child in {root/left,root/right,right/middle,right/deep,deep/lastleft,deep/lastright}
    \draw (\parent)--(\child);
  \end{scope}
  \node[anchor=south] at (7.2,.83) {rank cuts: active cells $A_j$};
  \foreach \j in {0,1,2,3} {
    \pgfmathsetmacro{\y}{.5-.75*(3-\j)}
    \node[anchor=east] at (4.9,\y-.175) {$j=\j$};
    \fill[black!4] (5.1,\y-.35) rectangle (9.5,\y);
    \ifnum\j<3
      \fill[omniBlue!18] (5.1,\y-.35) rectangle (7.3,\y);
      \draw (7.3,\y-.35)--(7.3,\y);
      \node at (6.2,\y-.175) {$a$};
    \fi
    \draw (5.1,\y-.35) rectangle (9.5,\y);
    \ifnum\j<2
      \draw (8.4,\y-.35)--(8.4,\y);
      \node at (7.85,\y-.175) {$c$};
    \fi
    \ifcase\j
      \draw (8.95,\y-.35)--(8.95,\y);
      \node at (8.675,\y-.175) {$e$};
      \node at (9.225,\y-.175) {$f$};
    \or
      \node at (8.95,\y-.175) {$d$};
    \or
      \node at (8.4,\y-.175) {$b$};
    \or
      \node at (7.3,\y-.175) {$v_0$};
    \fi
  }
  \draw[->,line width=.7pt] (10,-2.1)--(10,.45)
    node[midway,right,rotate=90,anchor=south] {coarsen};
\end{tikzpicture}
\caption{At rank $1$, $\{e,f\}$ becomes $d$, giving $A_1=\{a,c,d\}$.
The shaded leaf $a$ persists until the root $v_0$ becomes active at rank $3$.}
\label{fig:rank-cuts}
\end{figure}

Multigrid prolongation evaluates coarse polynomials on contained fine cells.
We construct the meshes and this containment relation by pruning one fixed
canonical finest tree $\Tcal$, preserving its grouping.
A \emph{leaf family} contains all children of a parent when they are leaves.
Each sweep replaces every such family by its parent. Rank is the sweep at
which a subtree becomes one leaf: rank zero is finest and increasing rank
coarsens (\Cref{fig:rank-cuts}). It measures remaining tree height, so cells
at one rank may have different sizes and coordinate levels. The cuts may
have unbalanced or crossed traces. The two-sided equations of
\Cref{sec:conforming-spaces} provide conforming FE spaces.
For a node $v$, define
\begin{equation}
\label{eq:node-rank-recurrence}
  r(v)=
  \begin{cases}
    0,&v\in\leaves(\Tcal),\\
    1+\max_{w\in\ch(v)}r(w),&v\notin\leaves(\Tcal).
  \end{cases}
\end{equation}
A node is active from rank $r(v)$ until its parent activates.
With $r(\parop(v_0))=\infty$ for the root $v_0$, set
\begin{equation}
\label{eq:active-rank-set}
  A_j:=\set{v\in\nodes(\Tcal)\given r(v)\le j<r(\parop(v))}.
\end{equation}
Let $U_j$ be the descriptor obtained by pruning the finest tree at the rank cut
$A_j$, and let $N_j$ denote its number of descriptor records.

At nested rank cuts, cell incidence pairs each fine cell with its containing
coarse cell. After adaptation it instead records positive-volume intersections.
We store one list per cell, using row offsets and contiguous entries
(compressed sparse row, CSR; \Cref{sec:fe-transfer}).
\begin{theorem}[Canonical rank cuts and their cell incidence]
\label{thm:rank-characterization}
\label{thm:linear-transfer}
Every subtree pruning $U'$ of a canonical descriptor $U$ is canonical.
With $A_j$ from \eqref{eq:active-rank-set}, simultaneous leaf-family pruning gives
\begin{equation}
\label{eq:rank-cut-conclusions}
 \Lcal(U_j)=\{\kappa(v):v\in A_j\},\qquad
 \operatorname{Can}(U_j)=U_j,\qquad 0\le j\le r(v_0).
\end{equation}
For $j\le k$, each $v\in A_j$ has exactly one ancestor
$\anc_{j,k}(v)\in A_k$. Ranks and child/subtree-end offsets take
$\Theta(N)$ postorder work and storage for $N$ finest-tree nodes.
Given these offsets, emitting $U_j$ takes optimal $\Theta(N_j)$ work.
For $j<k$, the cell incidence
$\Ical_{j,k}=\{(v,\anc_{j,k}(v)):v\in A_j\}$ and its inverse CSR rows
in preorder leaf order take optimal $\Theta(|A_j|+|A_k|)$ work and storage.
\end{theorem}
\begin{proof}
\emph{Ranks.} Induction in \eqref{eq:node-rank-recurrence} gives the
activation interval $r(v)\le j<r(\parop(v))$. This identifies the active
cells, nesting and the unique ancestor.

\emph{Canonicality.} At a surviving internal node $(\kappa,\sigma)$, let
$\Lcal'$ be the pruned leaves. For $a\notin\sigma$, let $\eta$ be an original
spanning leaf and $\zeta$ its retained ancestor.
By \eqref{eq:maximal-mask-characterization} and \eqref{eq:containment},
\[
 \begin{cases}
  \ell_a(\zeta)>\ell_a(\kappa)\quad\forall\zeta\in\Lcal',&a\in\sigma,\\
  \ell_a(\kappa)\le\ell_a(\zeta)\le\ell_a(\eta)=\ell_a(\kappa)
  \quad\text{for some }\zeta\in\Lcal',\ \kappa\preceq\zeta\preceq\eta,&a\notin\sigma.
 \end{cases}
\]
Thus $s(\kappa,\Lcal')=\sigma$, proving canonicality by \Cref{thm:normal-form}.

\emph{Work and cell incidence.} Postorder computes ranks and offsets in
$\Theta(N)$ work. Skipping discarded subtrees visits exactly the $N_j$
emitted nodes. Carrying the active $A_k$ ancestor emits $\Ical_{j,k}$.
Each coarse leaf's fine descendants are consecutive, giving its inverse
CSR row. At most $2|A_j|-1$ nodes are visited (\Cref{sec:cell-ordering}),
so both CSR directions cost $\Theta(|A_j|+|A_k|)$.
\end{proof}

Canonical cuts need no fresh reconstruction. These bounds use preorder leaf
order. Lexicographic ordering adds the cost in \Cref{sec:computational-model}.

\subsubsection*{Hierarchy output size}
Separately stored cuts repeat cells that survive several ranks.
Let $R$ be the largest macro-root rank and $M=|\mathcal M|$.
The macro-wise cut $\Fcal_j$ has $L_j$ leaves and $N_j$ descriptor records.
Quadratic output occurs when the coordinate-level bound grows: a binary tree
with $h$ internal nodes on one path has $N=2h+1$, $N_j=2(h-j)+1$,
\begin{equation}
\label{eq:spine-hierarchy-size}
 \sum_{j=0}^{h}N_j=(h+1)^2=\Theta(N^2).
\end{equation}
Selecting fewer cuts avoids this repeated output. The heat experiments
start at rank zero and repeatedly select the first coarser cut with at most
half as many leaves while $L_j>2M$, then append the root cut if absent
(\Cref{sec:marking-comparison}). For selected ranks $J$, the geometric decrease gives
\[
 \sum_{j\in J\setminus\{R\}}L_j\le2L_0,\qquad
 \sum_{j\in J}N_j\le4L_0.
\]
The supplement's \emph{Hierarchy counts} gives the exact record count and
an $O(N+R)$ prefix-sum construction of all leaf counts and selected ranks.
Computing ranks once and emitting selected cuts and cell-incidence CSR arrays
costs $O_{d,w}(N+\sum_{j\in J}N_j)$ work
(\Cref{thm:rank-characterization,sec:computational-model}).

\section{Face closures and complexity}
\label{sec:faces}
\label{sec:traces}

An overlay records opposite trace intersections without changing volume cells.
\emph{Closure} refines cells: balance limits tangential width ratios to two,
while optional trace comparability makes overlapping traces nested for
one-sided interpolation. Two-sided equations also handle crossings
(\Cref{sec:conforming-spaces}). We bound closure growth, then maintain
contacts and account for the complete mesh cycle.

\subsection{Trace overlays}
\label{sec:trace-overlays}

The two macro faces normal to axis $n$ are the coordinate planes
$f_{n,0}=\{\xi_n=0\}$ and $f_{n,1}=\{\xi_n=1\}$ of the reference hypercube.
A leaf touches the lower ($\varepsilon=0$) or upper ($\varepsilon=1$) face exactly when
\begin{equation}
\label{eq:face-incidence}
 i_n(\kappa)=\varepsilon(2^{\ell_n(\kappa)}-1).
\end{equation}
The leaf's \emph{trace} $F$ is its boundary portion on this face.
Its tangential key $\tr_f\kappa=((\ell_a,i_a))_{a\ne n}$ defines a box in
$[0,1]^{d-1}$ with level vector $\lambda(F)=(\lambda_j(F))_{j=1}^{d-1}$. Thus
\[
 \Pcal_f(\Tcal)=\{Q_{\tr_f\kappa(v)}:
 v\in\leaves(\Tcal),\ i_n(\kappa(v))=\varepsilon(2^{\ell_n(\kappa(v))}-1)\}.
\]
These half-open boxes partition the face and inherit an admissible descriptor
by retaining touching children, deleting normal splits and suppressing empty roots.
A leaf face $f=(n,\varepsilon)$ lies at $\xi_n=(i_n+\varepsilon)2^{-\ell_n}$.
Its trace uses macro coordinates and levels without rescaling.
In $d=1$ a trace is a point. A \emph{face contact} is a positive
$(d-1)$-measure contact within a macro or across a prescribed macro-face gluing.

A macro-face gluing \eqref{eq:signed-permutation} permutes coordinates and flips dyadic indices by
\begin{equation}
\label{eq:dyadic-flip}
  (\ell,i)\longmapsto (\ell,2^\ell-1-i),
\end{equation}
with inverse permutation and permuted flips for inverse gluing. We apply gluings
to keys, then take half-open representatives: literal endpoint reflections
differ only on sets of measure zero.

The gluing transports the second trace partition into the first's
reference-face chart.
Write $F_1\pitchfork F_2$ for overlapping relative interiors, excluding
lower-dimensional contacts. Their \emph{trace overlay} is
\begin{equation}
\label{eq:trace-overlay}
  \Pcal_1\wedge \Phi\Pcal_2
  :=\set{F_1\cap F_2\given F_1\in\Pcal_1,
  \ F_2\in\Phi\Pcal_2,
  \ F_1\pitchfork F_2}.
\end{equation}

\begin{figure}[!htbp]
\centering
\begin{tikzpicture}[x=1cm,y=0.85cm,>=Stealth,line join=round,
    cell/.style={draw,line width=.55pt},
    every node/.style={font=\small}]
  \begin{scope}
    \path[cell,fill=black!4] (0,1.6)--(1.6,1.6)--(1.05,1.95)--(-.55,1.95)--cycle;
    \path[cell,fill=omniBlue!18] (0,0)--(-.55,.35)--(-.55,1.95)--(0,1.6)--cycle;
    \path[cell,fill=omniBlue!18] (0,1.6)--(.8,1.6)--(.25,1.95)--(-.55,1.95)--cycle;
    \fill[omniBlue!18] (0,0) rectangle (.8,1.6);
    \draw[cell] (0,0) rectangle (1.6,1.6);
    \draw[cell,draw=omniBlue!85!black] (.8,0)--(.8,1.6);
    \node at (.4,.8) {$F$};
    \node[font=\footnotesize] at (.6,2.3) {side 1: vertical strips};
    \node[font=\footnotesize] at (.8,-.3) {$\lambda(F)=(1,0)$};
  \end{scope}
  \node at (2.25,.8) {$\wedge$};
  \begin{scope}[shift={(3.25,0)}]
    \path[cell,fill=black!4] (0,1.6)--(1.6,1.6)--(1.05,1.95)--(-.55,1.95)--cycle;
    \path[cell,fill=black!4] (0,0)--(-.55,.35)--(-.55,1.95)--(0,1.6)--cycle;
    \path[cell,fill=omniOrange!20] (0,0)--(-.55,.35)--(-.55,1.15)--(0,.8)--cycle;
    \draw[cell,draw=omniOrange!85!black] (0,.8)--(-.55,1.15);
    \fill[omniOrange!20] (0,0) rectangle (1.6,.8);
    \draw[cell] (0,0) rectangle (1.6,1.6);
    \draw[cell,draw=omniOrange!85!black] (0,.8)--(1.6,.8);
    \node at (.8,.4) {$G$};
    \node[font=\footnotesize] at (.6,2.3) {side 2: horizontal strips};
    \node[font=\footnotesize] at (.8,-.3) {$\lambda(G)=(0,1)$};
  \end{scope}
  \draw[->,line width=.7pt] (5.2,.8)--(6.2,.8);
  \begin{scope}[shift={(6.55,0)}]
    \fill[omniBlue!12] (0,.9) rectangle (.9,1.8);
    \fill[omniOrange!14] (.9,0) rectangle (1.8,.9);
    \fill[black!22] (0,0) rectangle (.9,.9);
    \draw[cell] (0,0) rectangle (1.8,1.8);
    \draw[cell,draw=omniBlue!85!black] (.9,0)--(.9,1.8);
    \draw[cell,draw=omniOrange!85!black] (0,.9)--(1.8,.9);
    \node[font=\footnotesize] at (.45,.45) {$F\cap G$};
    \node[font=\footnotesize] at (.9,2.3) {shared-face overlay};
    \node[font=\footnotesize] at (.9,-.3) {balanced, but crossed};
    \draw[->] (2.05,.15)--(2.7,.15) node[right] {$\zeta_1$};
    \draw[->] (2.05,.15)--(2.05,.95) node[above] {$\zeta_2$};
  \end{scope}
\end{tikzpicture}
\caption{Balanced, crossed traces in common tangential coordinates
$(\zeta_1,\zeta_2)$, with level difference $(1,-1)$. The overlay records
intersections without refining either volume mesh.}
\label{fig:trace-overlay}
\end{figure}
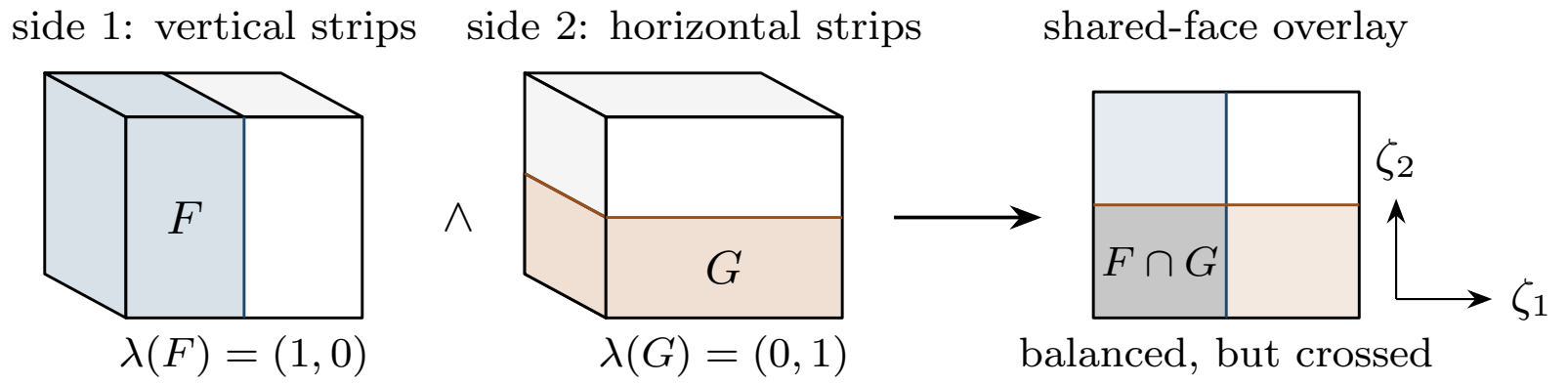

\begin{lemma}[Common refinement by trace intersections]
\label{thm:trace-overlay}
For two admissible trace partitions, the positive-measure intersections
$H=F\cap G$, $F\in\Pcal_1$, $G\in\Phi\Pcal_2$, form the least common
dyadic refinement \eqref{eq:trace-overlay}, which is admissible.
Each $H$ determines its factor pair $(F,G)$ uniquely.
The factors need not be nested for $d-1>1$.
\end{lemma}
\begin{proof}
\emph{Partition and leastness.} Dyadic nesting gives cells; the two
half-open partitions give disjoint coverage and unique factors.
For any common-refinement cell $C$,
$C\subseteq F,\ C\subseteq G\Longrightarrow C\subseteq F\cap G$.

\emph{Admissibility.} Restrict the second tree to each leaf $F$ of the first.
At a visited cell $B$, put $C=F\cap B$. Retain split axes whose midpoint
crosses $C$; on the others select the half containing $C$ and suppress empty
masks. Dyadic nesting makes retained cuts midpoints of $C$.
The leaves are exactly $F\cap G$. Grafting at $F$ gives the overlay
(\eqref{eq:child-partition}). Signed permutations preserve this construction;
\Cref{fig:trace-overlay} gives nonnested factors.
\end{proof}

Crossed traces of $n$ and $m$ strips give $nm$ overlay cells; the
shared-tree bound of \Cref{thm:linear-transfer} does not apply.

\subsection{Face balance}
\label{sec:balance}

A face-balanced forest has tangential width ratios at most two in common
face charts, while allowing arbitrary cell aspect ratios.
Normal splits preserve traces and unglued boundaries impose no condition.

Write $\iota_h(j)$ for the cell axis corresponding to tangential axis $j$ on face $h$.
We compare faces $f$ of $K$ and $f'$ of $L$ in chart $f$: the gluing $\Phi:f'\to f$
has permutation $\pi$ (reciprocally $\pi^{-1}$). The cell axes are
\begin{equation}
\label{eq:comparison-axis-maps}
  \jmath_F(j):=\iota_f(j),\qquad
  \jmath_G(j):=\iota_{f'}(\pi(j)).
\end{equation}
Flips affect positions, not axes. Internal contacts use $\pi=\mathrm{id}$.

\begin{definition}[Face balance]
\label{def:balance}
A forest is face-balanced if, for every face contact between transported
active traces $F,G$,
\begin{equation}
\label{eq:balance-condition}
  |\lambda_j(F)-\lambda_j(G)|\le 1,
  \qquad j=1,\ldots,d-1 .
\end{equation}
\end{definition}

For a face contact between traces $F$ of $K$ and $G$ of $L$, the balance
split demands are the coarser-side axes whose level gap exceeds one:
\begin{equation}
\label{eq:balance-demand}
  \begin{aligned}
  D_{\rm bal}(K;F,G)&:=\set{\jmath_F(j)\given \lambda_j(G)>\lambda_j(F)+1},\\
  D_{\rm bal}(L;G,F)&:=\set{\jmath_G(j)\given \lambda_j(F)>\lambda_j(G)+1}.
  \end{aligned}
\end{equation}
Algorithm~\ref{alg:balance-closure} unions all demands before splitting.
\Cref{sec:interface-output-complexity} implements these same synchronous
sweeps by retaining matched contacts.

\begin{algorithm}[htbp]
\caption{Balance by synchronous sweeps}
\label{alg:balance-closure}
{\small
\begin{algorithmic}[1]
\Require input forest \(\Fcal^0\)
\State \(\Fcal\gets\Fcal^0\)
\While{true}
  \State \(\sigma_K\gets\emptyset\) for every active cell \(K\)
  \For{each face contact between traces \(F,G\) of cells \(K,L\), counted once}
    \State \(\sigma_K\gets\sigma_K\cup D_{\rm bal}(K;F,G)\)
    \State \(\sigma_L\gets\sigma_L\cup D_{\rm bal}(L;G,F)\)
  \EndFor
  \State \textbf{if} all masks are empty \textbf{then return} the canonical descriptor of $\Fcal$
  \State simultaneously replace marked \(K\) by \(\ch_{\sigma_K}(K)\)
\EndWhile
\end{algorithmic}
}
\end{algorithm}

Write \(\Fcal\le_{\rm ref}\Hcal\) when \(\Hcal\) refines \(\Fcal\), that is,
when every active cell of \(\Hcal\) lies in an active cell of \(\Fcal\).

A demand $a=\jmath_F(j)$ in \eqref{eq:balance-demand} forces an $a$-split
on every branch of $K$ meeting $\relint(F\cap G)$ in any balanced refinement.
Away from trace boundaries, the branch cell $K'$ has trace $F'\subseteq F$
opposite $G'\subseteq G$, and balance gives
\begin{equation}
\label{eq:forced-axis-level}
 \ell_a(K')=\lambda_j(F')\ge\lambda_j(G')-1
 \ge\lambda_j(G)-1>\lambda_j(F)=\ell_a(K).
\end{equation}
Algorithm~\ref{alg:balance-closure} applies this branchwise demand to the whole
owning cell $K$.

\begin{theorem}[Balance closure under local level ceilings]
\label{thm:balance-closure}
Let the input \(\Fcal^0\) refine a reference forest \(\mathcal R\) on the same macro complex.
Each reference leaf \(K\) has integer ceilings
\(\ell_a(K)\le U_a(K)\le L_{\max}\), bounding all its input descendants.
Across every reference contact, transported tangential ceilings differ by
at most one. Then \Cref{alg:balance-closure} terminates with a deterministic,
face-balanced output $\Fcal^*$. Every intermediate descendant \(K'\subseteq K\) satisfies
\(\ell_a(K')\le U_a(K)\). The output leaf count obeys
\begin{equation}
\label{eq:closure-leaf-budget}
 \card{\Pcal(\Fcal^*)}\le
 \sum_{K\in\Pcal(\mathcal R)}2^{\sum_{a=1}^d(U_a(K)-\ell_a(K))}.
\end{equation}
In particular, closure never exceeds the maximal input coordinate level
\(L^0\). If \(\mathcal R\) is face-balanced and input descendants gain at most
\(r\) levels per axis, closure preserves this bound and has at most
\(2^{dr}\card{\Pcal(\mathcal R)}\) leaves.
\end{theorem}

\begin{proof}
Induct on sweeps. Let $U,V$ be the transported ceilings on a contact.
A demanded split satisfies
\[
 \lambda_j(F)+1\le\lambda_j(G)-1\le V-1\le U.
\]
Inherited reference contacts, including self-gluings, use compatible
transported ceilings. New internal contacts within one reference leaf have
$U=V$. Thus all splits preserve the ceilings.
Every current leaf is a union of ceiling-grid boxes. For
$N_r=|\Pcal(\Fcal^r)|$, this gives
\[
 N_r\le\sum_{K\in\Pcal(\mathcal R)}2^{\sum_a(U_a(K)-\ell_a(K))},
 \qquad N_{r+1}>N_r\quad\text{if any demand remains}.
\]
Thus closure terminates. Empty demands give balance and mask unions give determinism.
The choices $U_a(K)=L^0$ on $\mathcal R=\Fcal^0$, and
$U_a(K)=\min\{L_{\max},\ell_a(K)+r\}$ on balanced $\mathcal R$,
give the final bounds.
\end{proof}

\subsubsection*{Constructing the ceilings.}
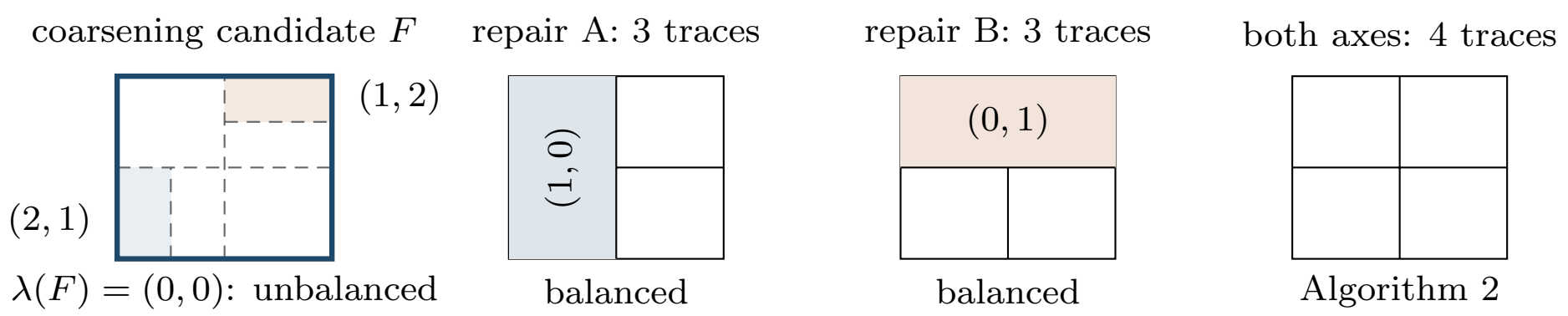
\begin{figure}[!htbp]
\centering
\begin{tikzpicture}[x=1cm,y=.85cm,>=Stealth,every node/.style={font=\footnotesize}]
  \foreach \x/\title in {0/{coarsening candidate $F$},3.1/{repair A: 3 traces},6.2/{repair B: 3 traces},9.3/{both axes: 4 traces}} {
    \node at (\x+.85,2.1) {\title};
    \draw (\x,0) rectangle (\x+1.7,1.7);
  }
  \fill[omniBlue!10] (0,0) rectangle (.425,.85);
  \fill[omniOrange!14] (.85,1.275) rectangle (1.7,1.7);
  \draw[densely dashed,black!60] (.85,0)--(.85,1.7) (0,.85)--(1.7,.85)
    (.425,0)--(.425,.85) (.85,1.275)--(1.7,1.275);
  \draw[omniBlue!85!black,line width=1.2pt] (0,0) rectangle (1.7,1.7);
  \node[anchor=east] at (-.08,.35) {$(2,1)$};
  \node[anchor=west] at (1.78,1.5) {$(1,2)$};
  \node at (.85,-.3) {$\lambda(F)=(0,0)$: unbalanced};
  \fill[omniBlue!15] (3.1,0) rectangle (3.95,1.7);
  \draw (3.95,0)--(3.95,1.7) (3.95,.85)--(4.8,.85);
  \node[rotate=90] at (3.525,.85) {$(1,0)$};
  \node at (3.95,-.3) {balanced};
  \fill[omniOrange!17] (6.2,.85) rectangle (7.9,1.7);
  \draw (6.2,.85)--(7.9,.85) (7.05,0)--(7.05,.85);
  \node at (7.05,1.275) {$(0,1)$};
  \node at (7.05,-.3) {balanced};
  \draw (10.15,0)--(10.15,1.7) (9.3,.85)--(11,.85);
  \node at (10.15,-.3) {Algorithm~\ref{alg:balance-closure}};
\end{tikzpicture}
\caption{Coarse face $F$ against fixed opposite traces (dashed).
The three-trace repairs A and B are balanced and incomparable. Their
overlapping shaded strips force the unbalanced whole $F$ into every common
coarsening, so no balanced repair is coarser than every other repair.
Algorithm~\ref{alg:balance-closure} uses four traces, so it need not minimize
their number.}
\label{fig:nonminimal-balance}
\end{figure}

Ceilings must dominate input levels and vary by at most one across each
tangential contact. The reference forest $\mathcal R$ gives an undirected graph $G$ with vertices
$v=(K,a)$ and edges joining transported tangential pairs across contacts.
Let $b_v$ be the largest $a$-level of an input descendant of $K$,
$C(v)$ its connected component and $\operatorname{dist}_G$ shortest-path distance. Then
\begin{equation}
\label{eq:least-compatible-ceilings}
 U_v^*=\max_{z\in C(v)}\bigl(b_z-\operatorname{dist}_G(v,z)\bigr),
 \qquad U_v^{\rm cmp}=\max_{z\in C(v)}b_z.
\end{equation}
For adjacent $v,w$ and $z\in C(v)$, the choice $z=v$, the triangle inequality
and compatibility of any majorant $U\ge b$ give, respectively,
\[
 U_v^*\ge b_v,\qquad U_v^*\le U_w^*+1,\qquad
 U_v\ge U_z-\operatorname{dist}_G(v,z)
      \ge b_z-\operatorname{dist}_G(v,z).
\]
Taking the maximum over $z$ in the last inequality proves $U\ge U^*$.
Component maxima give the least majorant with equal neighboring ceilings.
$U^*$ supplies balance ceilings, and $U^{\rm cmp}$ supplies ceilings for
balance followed by comparability
(\Cref{thm:trace-comparability-closure}). They minimize the corresponding
grid budget \eqref{eq:closure-leaf-budget} for this reference forest.
Ceiling minimality does not imply a least balanced mesh repair: for $d\ge3$,
whole-cell demands can give incomparable repairs
(\Cref{fig:nonminimal-balance}).

For one-level refinement of balanced input, the least ceilings identify
where balance can propagate. On the fixed input cells of cycle $t$, draw $(K,a)\to(L,b)$
across a contact when $a,b$ are corresponding tangential axes and
$\ell_b(L)=\ell_a(K)-1$. Let $k_{t,K}$ count the axes $a$ of $K$ for
which $(K,a)$ is reachable from any requested pair, including the requests
themselves. Reachability bounds the possible split support.

\subsubsection*{Propagation sheets.}
For $K=(m,\kappa)$, write $I_a(K)=I_{\ell_a(K),i_a(K)}$.
The \emph{sheet} $\Sigma(K,a)$ is the smallest set containing $(m,a,I_a(K))$
and closed under transport across macro faces whose normal differs from the
transported axis, using \eqref{eq:signed-permutation}.
Each triple specifies a coordinate slab in one macro.
At most $2d$ triples occur per macro: the interval is unchanged or reflected.
Let $\sigma_t(K,a)$ count the old cell--axis pairs $(L,b)$ whose $I_b(L)$
contains a sheet interval in that macro and axis, counting each pair once.
Transport can return to another axis of the same cell.

\Needspace{8\baselineskip}
\begin{theorem}[Growth controlled by requested bisections]
\label{thm:marked-growth}
Let $d\ge2$ and let $\Fcal^0$ be face-balanced. Each cycle $t<T$ requests
one-level refinement followed by balance alone. Let $N_t$ be the leaf count and
$\mathcal A_t$ the set of requested cell--axis pairs, with $\ell_a(K)$ measured
before refinement. In each cycle,
\begin{equation}
\label{eq:interaction-marked-growth}
 N_{t+1}-N_t\le
 \sum_{K\in\Pcal(\Fcal^t)}\bigl(2^{k_{t,K}}-1\bigr).
\end{equation}
This counts each reached coordinate once, including when propagation from
different requests overlaps. Uniform all-axis refinement has $k_{t,K}=d$
and attains equality, with increase $(2^d-1)N_t$.
An explicit worst-case bound uses
\begin{equation}
\label{eq:marking-depth-weight}
 D_d=2(d-1)2^{d-2},\qquad
 B_d=D_d-(d-1),\qquad
 H_d(\ell)=1+D_d\sum_{j=0}^{\ell-1}B_d^j.
\end{equation}
Here $D_d,B_d$ bound initial and subsequent branching of balance demands.
A request $(K,a)$ reaches at most
$\min\{H_d(\ell_a(K)),\sigma_t(K,a)\}$ cell--axis pairs.
Equation~\eqref{eq:interaction-marked-growth} gives
\begin{equation}
\label{eq:cumulative-marked-growth}
 N_T-N_0\le\frac{2^d-1}{d}
   \sum_{t<T}\ \sum_{(K,a)\in\mathcal A_t}\min\{H_d(\ell_a(K)),\sigma_t(K,a)\}.
\end{equation}
For $d=2$, at most two descending chains give $H_2(\ell)=1+2\ell$.
\end{theorem}

\begin{proof}
Fix a cycle and write $\ell_v=\ell_a(K)$ for $v=(K,a)$.
The post-request bounds are $b_v=\ell_v+\1_{\mathcal A_t}(v)$.
Input balance gives $\ell_z-\operatorname{dist}_G(v,z)\le\ell_v$.
Equality requires a path decreasing by one old level at every edge. Hence
\eqref{eq:least-compatible-ceilings} gives
\[
 \begin{aligned}
 U_v^*-\ell_v
 &=\max_{z\in C(v)}\bigl(\1_{\mathcal A_t}(z)+\ell_z-\ell_v
                       -\operatorname{dist}_G(v,z)\bigr)\\
 &=\1_{\{v\text{ reachable from }\mathcal A_t\}}.
 \end{aligned}
\]
By \Cref{thm:balance-closure}, $K$ has at most $2^{k_{t,K}}$ descendants.
Summing their increases gives \eqref{eq:interaction-marked-growth}.

The half-grid count in \Cref{app:propagation-branching} gives outdegree
$D_d$, reduced to $B_d$ after an incoming edge. A path from $(K,a)$ has at
most $\ell_a(K)$ edges, so its reachable set has size at most
$1+D_d\sum_{j=0}^{\ell_a(K)-1}B_d^j=H_d(\ell_a(K))$.

The propagated axis $b$ is tangential to its contact, so an edge crosses only
a face normal to another axis. Dyadic nesting gives, in the common face chart,
\[
 \ell_{b'}(L')=\ell_b(L)-1,\qquad I_{b'}(L')\supseteq I_b(L).
\]
Identity transport within macros and signed transport across macro faces
keep every reachable vertex on the sheet, giving the bound $\sigma_t(K,a)$.

Finally, $2^k-1\le(2^d-1)k/d$ for $0\le k\le d$, and the union bound gives
\[
 \sum_K k_{t,K}\le
 \sum_{(K,a)\in\mathcal A_t}\min\{H_d(\ell_a(K)),\sigma_t(K,a)\}.
\]
Summing \eqref{eq:interaction-marked-growth} over cycles proves
\eqref{eq:cumulative-marked-growth}.
\end{proof}

\begin{corollary}[Balance sweep count]
\label{cor:balance-sweeps}
Under the hypotheses of \Cref{thm:marked-growth}, the number $S_t$ of
nonterminal balance sweeps satisfies
\begin{equation}
\label{eq:balance-sweep-budget}
 S_t\le\max_{(K,a)\in\mathcal A_t}\ell_a(K),\qquad \max\emptyset=0.
\end{equation}
\end{corollary}
\begin{proof}
By the preceding proof, each old cell--axis pair gains at most one level.
New internal contacts have gaps at most one. Orient a demanded coordinate
from the finer current trace $F$ to the coarser $G$. With pre-cycle signed gap
$\delta\le1$ and cumulative increments $\varepsilon_F,\varepsilon_G\in\{0,1\}$,
\[
 \delta+\varepsilon_F-\varepsilon_G=2
 \quad\Longrightarrow\quad
 (\delta,\varepsilon_F,\varepsilon_G)=(1,1,0).
\]
After sweep one, the finer-side increment must occur in the preceding sweep.
Otherwise the parent contact already demanded the split. Writing $\ell_v$
for the old level of pair $v$, a demand in sweep $r$ therefore has a causal chain
\[
 v_0\in\mathcal A_t,\quad v_0\to\cdots\to v_r,\qquad
 0\le\ell_{v_r}=\ell_{v_0}-r.
\]
This proves \eqref{eq:balance-sweep-budget}.
\end{proof}

\subsubsection*{Sheets and depth.}
Unlike simplicial completion
\cites{stevensonCompletion2008,dieningGehringStornRefinement2026}, this whole-cell
rule has no depth-independent cumulative marking bound. In 2D, $2s-1$
bisections create $2s$ vertical strips; $s$ further bottom-cell $y$-bisections
in the central strip produce $s^2+2s$ leaves.
For $d\ge3$, a prepared balanced cube has $L+1$ dyadic $x$-layers.
Layer $h$ carries $2^{(d-2)h}$ transverse columns with bottom $y$-depth $L-h$.
One layer-zero bottom-cell request $(K,a)$ attains the sheet count:
\[
 \Delta N=\sigma_t(K,a)=\sum_{h=0}^L2^{(d-2)h},\qquad S=L.
\]
This excludes preparation and requires $L_{\max}\ge L+1$.
Both families use balance alone and allow crossed traces. The supplement's
\emph{Growth constructions} gives exact histories and normal-graded repairs:
$1+2^{d-1}$ added leaves for the layer family, and linear growth for the strip
history when extra normal levels are available. The whole-cell policy refines
only axes demanded by contact imbalances; \Cref{tab:large-cycles} reports its
closure factor $\gamma$. The gap between upper base $B_d$ and lower base
$2^{d-2}$ remains open. Logical sweeps give no communication-round lower bound
\cite{isaacBursteddeGhattasBalance2012}.

\subsection{Trace comparability for one-sided constraints}
\label{sec:trace-comparability}

Optional closure supplies nested traces for one-sided interpolation.
Two-sided equations handle crossings (\Cref{thm:conforming-kernel}).
For overlapping traces and face $f=(n,\varepsilon)$ of $K=(m,\kappa)$,
\begin{align}
\label{eq:trace-comparability}
  F\subseteq G
  &\quad\Longleftrightarrow\quad
  \lambda(F)\ge\lambda(G)\ \text{componentwise},\\
\label{eq:face-identity}
 \operatorname{id}(K,f)&=(m,(\ell_a)_a,(i_a)_a,n,\varepsilon).
\end{align}
A contact is \emph{crossed} if neither inclusion holds, and
\emph{trace-comparable} otherwise. Full identities distinguish self-gluings.
For a crossed face contact $e=(F,G)$ owned by $K,L$, use the axis maps
\eqref{eq:comparison-axis-maps} to define its missing-axis masks
\[
 M_K(F,G)=\{\jmath_F(j):\lambda_j(F)<\lambda_j(G)\},\qquad
 M_L(G,F)=\{\jmath_G(j):\lambda_j(G)<\lambda_j(F)\}.
\]
The trace with fewer missing axes is selected. Ties use the
lexicographically smaller full face identity.
$\Delta_{\rm cmp}(e,K)$ is its mask at the selected owner and empty elsewhere.
For the crossed contacts $\mathcal X_r$ of a balanced forest $\Fcal^r$, set
\begin{equation}
\label{eq:compatibility-iteration}
 \sigma_K^r=\bigcup_{e\in\mathcal X_r}\Delta_{\rm cmp}(e,K),\qquad
 \Fcal^{r+1}=\operatorname{Bal}(\operatorname{Split}(\Fcal^r,\sigma^r)).
\end{equation}
Split applies the masks simultaneously; Bal is
Algorithm~\ref{alg:balance-closure}. Iteration stops when $\mathcal X_r=\emptyset$.

The comparability bound uses equal transported ceilings: a missing-axis split
raises the selected coordinate to the opposite trace's level.
Matching reference traces supply these ceilings in the final case below.

\begin{theorem}[Trace-comparability closure under local level ceilings]
\label{thm:trace-comparability-closure}
Under \Cref{thm:balance-closure} with \emph{equal} transported ceilings and
face-balanced input, \eqref{eq:compatibility-iteration} terminates with a
balanced, trace-comparable output, preserving the ceilings and leaf budget.
The stated tie-breaking policy gives determinism. If reference traces match
after chart transport and input descendants gain at most $r$ levels per axis, the output has at most
$2^{dr}\card{\Pcal(\mathcal R)}$ leaves.
\end{theorem}

\begin{proof}
For a selected tangential coordinate with common transported ceiling $U$,
$\lambda_j(F)+1\le\lambda_j(G)\le U$.
Thus Split preserves the ceilings, as does Bal by \Cref{thm:balance-closure}.
The finite leaf budget and strict growth give termination. Mask unions and
full face identities give determinism.
Constant ceilings give the global bound. Matching reference traces permit
$U_a(K)=\min\{L_{\max},\ell_a(K)+r\}$.
\end{proof}

\subsection{Geometric bounds for balance and adaptation}
\label{sec:balance-complexity}

Level gaps bound overlap counts: subdivide a query box until dyadic nesting
places each subbox in one cell, then locate their centers.

\Needspace{11\baselineskip}
\begin{lemma}[Dyadic packing and probes]
\label{lem:final-overlap-count}
\label{lem:dyadic-probing}
Let a dyadic partition $\Pcal$ cover a dyadic $k$-box $K$, and suppose
$\ell_a(\eta)\le\ell_a(K)+r_a$ for every positive-measure intersection
$K\cap\eta$, with $r_a\in\Nzero$. Subdivide $K$ by $r_a$ levels per axis
to obtain $\mathcal Q$, with centers $c_Q$. Then
\begin{gather}
\label{eq:dyadic-probe-coverage}
 K\cap\eta=\dotcup_{Q\in\mathcal Q:\,Q\subseteq\eta}Q,\qquad
 \{\eta:|K\cap\eta|>0\}=\{\eta(c_Q):Q\in\mathcal Q\},\\
\label{eq:face-overlap-packing}
 \card{\{\eta:|K\cap\eta|>0\}}\le|\mathcal Q|=2^{\sum_a r_a},
\end{gather}
where $\eta(x)$ is the unique cell of $\Pcal$ containing $x$.
This bound is sharp. A face-balanced forest therefore has at most
$2^{d-1}$ opposite traces per leaf face and $O_d(N)$ contacts for $N$ nodes.
\end{lemma}
\begin{proof}
For $Q\in\mathcal Q$, dyadic nesting \eqref{eq:dyadic-dichotomy} and the level bound give
$Q\cap\eta\ne\emptyset\Longleftrightarrow Q\subseteq\eta
\Longleftrightarrow c_Q\in\eta$.
Since $\mathcal Q$ partitions $K$, this proves
\eqref{eq:dyadic-probe-coverage}--\eqref{eq:face-overlap-packing}.
Taking $\Pcal=\mathcal Q$ gives equality.
For balanced forests use $k=d-1$, $r_a=1$ on at most $2dN$ faces.
For $k=0$, one probe and counting measure give the same argument.
\end{proof}

\begin{proposition}[Input-to-final incidence and recorded origins]
\label{prop:adaptation-incidence}
\label{prop:source-transfer}
Let $\mathcal I(K)$ be the input active cells meeting a final active cell $K$
in positive volume after one admissible adaptation and further refinement.
For $L_*$ final active cells,
\begin{equation}
\label{eq:input-final-packing}
 |\mathcal I(K)|\le2^d,\qquad \sum_K|\mathcal I(K)|\le2^dL_*.
\end{equation}
Both bounds are sharp, without balance or a bound on further refinement depth.
Also $\ell_a(\eta)\le\ell_a(K)+1$ for $\eta\in\mathcal I(K)$ and every $a$.
With supplied input indices and immutable origin expressions, both incidence
CSR directions take $O_d(n+n_*)$ additional work and storage, where $n,n_*$
count input and final descriptor records.
\end{proposition}
\begin{proof}
\emph{Overlap bound.} By \Cref{lem:lineage-invariant},
$\ell_a(\eta)\le\ell_a(K)+1$ for $\eta\in\mathcal I(K)$, also after further
refinement. Thus \Cref{lem:dyadic-probing} with $r_a=1$ gives
$|\mathcal I(K)|\le2^d$. Summation gives \eqref{eq:input-final-packing};
coarsening all $2^d$ root children attains equality.

\emph{Recorded origins.} Initially $\mathcal O(K_i)=\{i\}$.
Refinement copies origins, and contraction uses the multiset union
$\mathcal O(K_0\cup K_1)=\mathcal O(K_0)\uplus\mathcal O(K_1)$.
Distinct contraction axes bound the expression depth by $d$
(\Cref{lem:lineage-invariant}). Induction gives
\[
 K\subseteq\bigcup_{i\in\mathcal O(K)}K_i,\qquad
 |\mathcal O(K)|\le2^d,\qquad
 \mathcal I(K)=\{K_i:i\in\mathcal O(K),\ |K\cap K_i|>0\}.
\]
Expansion, filtering and deduplication cost $O_d(1)$ per final row.
Counting and transposition cost
$O(n+L_*+\sum_K|\mathcal I(K)|)=O_d(n+n_*)$.
Contractions create $O_d(n)$ shared union records; refinement stores
$O(n_*)$ handles (\Cref{sec:cell-ordering}).
Sorting and bit costs remain separate.
\end{proof}

To apply \eqref{eq:face-overlap-packing} during closure, we must also control
temporary level gaps. Define
$g_0(\Fcal)=\max_{(F,G),j}|\lambda_j(F)-\lambda_j(G)|$, taking
$\max\emptyset=0$.

\Needspace{10\baselineskip}
\begin{proposition}[Level gaps during adaptation and closure]
\label{prop:one-step-imbalance}
\label{prop:transient-gap}
One admissible adaptation of a forest $\Fcal$ gives $\Fcal'$ with
$g_0(\Fcal')\le g_0(\Fcal)+2$. Face-balanced input gives $g_0(\Fcal')\le3$.
Starting from any forest $\Fcal^0$, apply
Algorithm~\ref{alg:balance-closure}, optionally followed by
\eqref{eq:compatibility-iteration}. Every intermediate forest $\Fcal^r$ satisfies
\[
 g_0(\Fcal^r)\le g:=\max\{g_0(\Fcal^0),2\},\qquad
 \text{neighbors per leaf face}\le2^{g(d-1)}.
\]
\end{proposition}
\begin{proof}
\emph{Adaptation.} Away from the old lower-dimensional skeleton, the input
neighborhoods of a new contact either share an old face (including
self-gluings) or lie in one input leaf.
For an old contact, let $\delta_j,\delta_j'$ be its old and new signed
gaps, with level changes $\Delta_j^F,\Delta_j^G$. Then
\[
 |\delta_j'|\le|\delta_j|+|\Delta_j^F|+|\Delta_j^G|
 \le g_0(\Fcal)+2
\]
by \Cref{lem:lineage-invariant}. Inside one input leaf, initially refined
axes gain one level everywhere, other axes change by $-1$ or $0$. Hence
\[
 \delta_j=0,\quad
 (\Delta_j^F,\Delta_j^G)\in\{(1,1)\}\cup\{-1,0\}^2
 \quad\Longrightarrow\quad |\delta_j'|\le1.
\]

\emph{Closure.} Sibling contacts have zero gap. All others restrict parent
face overlaps by \eqref{eq:child-partition}, with self-gluings distinguished
by \eqref{eq:face-identity}.
For $\delta_j=\lambda_j(F)-\lambda_j(G)$,
$\delta_j'=\delta_j+\varepsilon_j^F-\varepsilon_j^G$,
$\varepsilon_j^F,\varepsilon_j^G\in\{0,1\}$.
The demands \eqref{eq:balance-demand} give
$\delta_j>1\Rightarrow\varepsilon_j^G=1$ and
$\delta_j<-1\Rightarrow\varepsilon_j^F=1$, so
\begin{equation}
\label{eq:transient-gap-recurrence}
 |\delta_j'|\le
 \begin{cases}|\delta_j|,&|\delta_j|>1,\\2,&|\delta_j|\le1\end{cases}
 \le\max\{|\delta_j|,2\}.
\end{equation}
Trace-comparability sweeps start balanced and therefore create gaps at most
two; balance preserves this bound. Induction gives $g_0(\Fcal^r)\le g$.
\Cref{lem:dyadic-probing} with $k=d-1$, $r_a=g$ gives the degree bound.
\end{proof}

\subsection{Constructing cell incidence and face contacts}
\label{sec:interface-output-complexity}

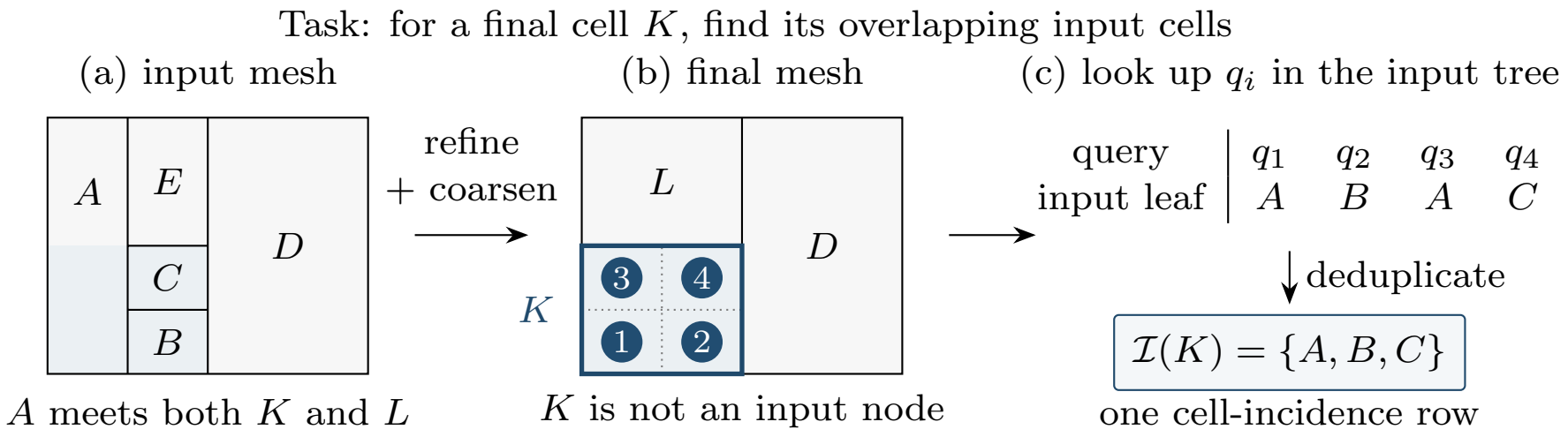
\begin{figure}[!htbp]
\centering
\begin{tikzpicture}[x=1cm,y=.8cm,>=Stealth,
  every node/.style={font=\footnotesize},
  probe/.style={circle,fill=omniBlue!90!black,text=white,inner sep=1pt,
    minimum size=3mm,font=\scriptsize}]
  \node at (5.3,3.25) {Task: for a final cell $K$, find its overlapping input cells};
  \node at (1.2,2.8) {(a) input mesh};
  \fill[black!3] (0,0) rectangle (2.4,2.4);
  \fill[omniBlue!9] (0,0) rectangle (1.2,1.2);
  \draw (0,0) rectangle (2.4,2.4);
  \draw (1.2,0)--(1.2,2.4) (.6,0)--(.6,2.4)
    (.6,1.2)--(1.2,1.2) (.6,.6)--(1.2,.6);
  \node at (.3,1.7) {$A$};
  \node at (.9,.3) {$B$}; \node at (.9,.9) {$C$};
  \node at (.9,1.8) {$E$}; \node at (1.8,1.2) {$D$};
  \node at (1.2,-.35) {$A$ meets both $K$ and $L$};
  \draw[->] (2.75,1.3)--(3.6,1.3);
  \node[align=center] at (3.175,1.95) {refine\\+ coarsen};
  \node at (5.2,2.8) {(b) final mesh};
  \begin{scope}[shift={(4,0)}]
    \fill[black!3] (0,0) rectangle (2.4,2.4);
    \fill[omniBlue!9] (0,0) rectangle (1.2,1.2);
    \draw (0,0) rectangle (2.4,2.4);
    \draw (1.2,0)--(1.2,2.4) (0,1.2)--(1.2,1.2);
    \draw[densely dotted,black!50] (.6,0)--(.6,1.2) (0,.6)--(1.2,.6);
    \draw[omniBlue!85!black,line width=1pt] (0,0) rectangle (1.2,1.2);
    \node[anchor=east,text=omniBlue!85!black] at (-.08,.6) {$K$};
    \node at (.6,1.8) {$L$}; \node at (1.8,1.2) {$D$};
    \foreach \x/\y/\n in {.3/.3/1,.9/.3/2,.3/.9/3,.9/.9/4} {
      \node[probe] at (\x,\y) {$\n$};
    }
  \end{scope}
  \node at (5.2,-.35) {$K$ is not an input node};
  \draw[->] (6.75,1.3)--(7.4,1.3);
  \node at (9.3,2.8) {(c) look up $q_i$ in the input tree};
  \node at (9.3,1.85) {$\begin{array}{c|cccc}
    \text{query}&q_1&q_2&q_3&q_4\\
    \text{input leaf}&A&B&A&C
    \end{array}$};
  \draw[->] (9.3,1.15)--(9.3,.65) node[midway,right] {deduplicate};
  \node[draw=omniBlue!85!black,fill=omniBlue!5,rounded corners=1pt,inner sep=4pt]
    at (9.3,.2) {$\mathcal I(K)=\{A,B,C\}$};
  \node at (9.3,-.35) {one cell-incidence row};
\end{tikzpicture}
\caption{One mixed adaptation splits $A$ horizontally and coarsens to
$K=A_{\rm low}\cup B\cup C$, $L=A_{\rm high}\cup E$. $D$ stays unchanged.
To recover $\mathcal I(K)$, queries in the input tree return $A,B,A,C$.
Dotted lines place queries, without further refinement.
The input tangential gap is two: cell-incidence probing does not require
face balance (\Cref{lem:dyadic-probing,prop:adaptation-incidence}).}
\label{fig:dyadic-probes}
\end{figure}

A final cell need not be an input-tree node (\Cref{fig:dyadic-probes}).
A \emph{probe} follows masks and branch bits to an input leaf.
The $2^d$ half-box centers recover all input overlaps
(\Cref{prop:adaptation-incidence,lem:dyadic-probing}). Face contacts instead
use recursive routing; matched records persist through closure.

The supplement's \emph{Gap-adaptive face probes} gives an alternative probe construction
and its exact coverage certificate when $g$ is unknown.

For a forest with $N$ descriptor records and branch length $D\le dL_{\max}$,
recursive routing emits $N_{\rm c}^0$ contacts with work
\begin{equation}
\label{eq:recursive-contact-work}
 T_{\rm match}\le C_d\!\left(N+(2D+1)N_{\rm c}^0\right).
\end{equation}
The supplement's \emph{Recursive routing bound} proves this estimate for the
uncached constructor used in the numerical experiments.

\subsubsection*{Retained-contact sweeps.}
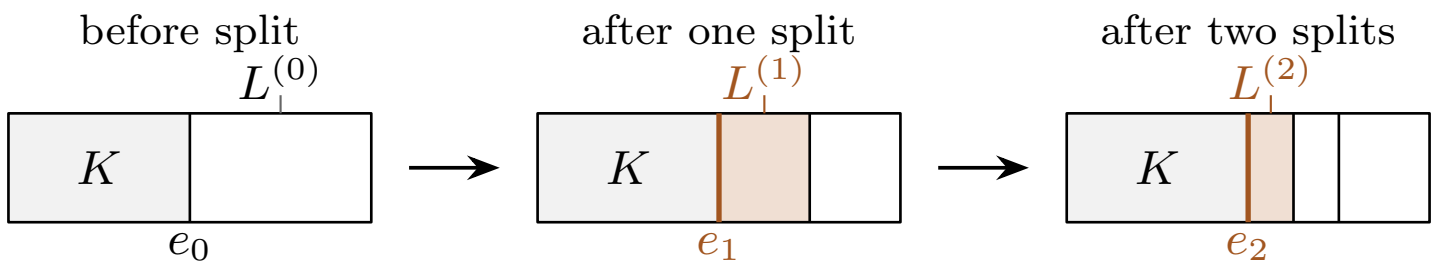
\begin{figure}[!htbp]
\centering
\begin{tikzpicture}[x=1cm,y=.6cm,>=Stealth,line join=round,
    every node/.style={font=\small},
    cell/.style={draw=black,line width=.5pt}]
  \foreach \epoch/\x/\paneltitle in {0/0/{before split},1/3.5/{after one split},2/7/{after two splits}} {
    \begin{scope}[shift={(\x,0)}]
      \fill[black!5] (0,0) rectangle (1.2,1.2);
      \ifnum\epoch>0
        \pgfmathsetmacro{\neighborend}{1.2+1.2/pow(2,\epoch)}
        \fill[omniOrange!20] (1.2,0) rectangle (\neighborend,1.2);
      \fi
      \draw[cell] (0,0) rectangle (2.4,1.2);
      \draw[cell] (1.2,0)--(1.2,1.2);
      \ifnum\epoch>0
        \draw[cell] (1.8,0)--(1.8,1.2);
        \draw[omniOrange!90!black,line width=1pt] (1.2,0)--(1.2,1.2);
      \fi
      \ifnum\epoch>1
        \draw[cell] (1.5,0)--(1.5,1.2);
      \fi
      \node at (.6,.6) {$K$};
      \pgfmathsetmacro{\neighborcenter}{1.2+.6/pow(2,\epoch)}
      \ifnum\epoch=0
        \draw[black!60] (\neighborcenter,1.2)--(\neighborcenter,1.4);
        \node at (\neighborcenter,1.6) {$L^{(\epoch)}$};
        \node at (1.2,-.27) {$e_\epoch$};
      \else
        \draw[omniOrange!90!black] (\neighborcenter,1.2)--(\neighborcenter,1.4);
        \node[text=omniOrange!90!black] at (\neighborcenter,1.6) {$L^{(\epoch)}$};
        \node[text=omniOrange!90!black] at (1.2,-.27) {$e_\epoch$};
      \fi
      \node[font=\footnotesize] at (1.2,2.08) {\paneltitle};
    \end{scope}
  }
  \draw[->,line width=.7pt] (2.65,.6)--(3.25,.6);
  \draw[->,line width=.7pt] (6.15,.6)--(6.75,.6);
\end{tikzpicture}
\caption{Normal refinement preserves $K$'s face but
replaces its neighbor $L^{(r)}$, retiring $e_r$ and creating $e_{r+1}$.
The replacement record is charged to the orange newborn neighbor
(\Cref{lem:contact-genealogy}).}
\label{fig:contact-birth-charge}
\end{figure}

A contact record binds immutable endpoint versions. It is \emph{live} until
an endpoint splits, then \emph{retired} (\Cref{fig:contact-birth-charge}).
A passed predicate remains true while the record is live. For
$p\in\{\mathrm{bal},\mathrm{cmp}\}$, let $\mathcal E_r^p$ contain the untested
live records at sweep start. Here $\Delta_{\rm bal}(e,K)$ unites the demands
\eqref{eq:balance-demand} at $K$, including both faces of a self-gluing.
$\Delta_{\rm cmp}$ is defined in \eqref{eq:compatibility-iteration}:
\begin{equation}
\label{eq:retained-sweep}
 \sigma_K^r=\bigcup_{e\in \mathcal E_r^p}\Delta_p(e,K),\qquad
 \Fcal^{r+1}=\operatorname{Split}(\Fcal^r,\sigma^r).
\end{equation}
Only touched masks are accumulated and cleared. After simultaneous commitment,
each split cell's old incidence list is visited. The first visit retires each
parent record and emits its overlapping child-face pairs once. Sibling
contacts and replacements are queued for the next sweep. Balance completes
before each comparability sweep. State persists until final canonical reconstruction.

\begin{lemma}[Contact genealogy]
\label{lem:contact-genealogy}
In \eqref{eq:retained-sweep}, assume tangential gaps stay bounded by $g$.
Each contact record has a unique birth certificate and at most one test per
predicate. Each cell birth receives at most $C_{d,g}$ contact charges.
\end{lemma}
\begin{proof}
Initial records use the fixed orientation of \eqref{eq:face-identity}.
A new record is certified by a sibling split or a retired parent record,
together with its overlapping child-face pair.
The child partition \eqref{eq:child-partition} makes the parent unique.
Retirement emits each pair once, even when both endpoints split.
Full face/gluing identities distinguish self-gluings.

Give initial versions birth sweep zero and order cell versions by birth
sweep and full cell key. For $e=(v,w)$, set
$\operatorname{charge}(e)=\max\{v,w\}$ in this order, including loops.
All records charged to $v$ exist at its birth-sweep end: later replacements
have a newborn endpoint $w>v$ (\Cref{fig:contact-birth-charge}).
Thus \Cref{lem:dyadic-probing} on the $2d$ faces gives
\begin{equation}
\label{eq:contact-charge}
 \card{\{e:\operatorname{charge}(e)=v\}}
 \le\deg_{\operatorname{birth}(v)}(v)\le C_{d,g}:=2d\,2^{g(d-1)}.
\end{equation}
Immutable endpoint versions keep each predicate unchanged until retirement.
Passing records are not retested. Failure forces replacement, and retired
records are not tested.
\end{proof}

Let $n_1,n_*$ count descriptor records at the start and end of closure,
and $T_{\rm upd}$ its contact updates and flat emission, excluding fresh matching.
\begin{corollary}[Retained-contact closure]
\label{cor:maintained-incidence}
At fixed dimension and key width,
the rule \eqref{eq:retained-sweep} implements Algorithm~\ref{alg:balance-closure}
and \eqref{eq:compatibility-iteration}
in $O_{d,w}(n_*)$ update/emission work and workspace.
Retirement and replacement visits are charged separately.
Closure only refines, so $n_1=O(n_*)$.
\end{corollary}
\begin{proof}
Set $g=\max(g_0,2)$ for initial gap $g_0$.
A passed predicate stays valid until an endpoint changes, so queuing new
contacts reproduces the full synchronous masks.

Let $N_{\rm b},N_{\rm c},N_{\rm t}$ count births, contact records and tests.
$N_{\rm c}^{\rm new}$ excludes initial contacts. For split arities $c_s\ge2$
and initial/final leaf counts $L_1,L_*$,
\begin{equation}
\label{eq:incidence-birth-budget}
\begin{aligned}
 N_{\rm b}
 &=\sum_s c_s\le2\sum_s(c_s-1)=2(L_*-L_1),\\
 N_{\rm c}^{\rm new}&\le C_{d,g}N_{\rm b},\qquad
 N_{\rm c}\le C_{d,g}(n_1+N_{\rm b}),\qquad N_{\rm t}\le2N_{\rm c}.
\end{aligned}
\end{equation}
By \Cref{lem:contact-genealogy}, \Cref{thm:marked-growth} also controls
$N_{\rm c}^{\rm new}$ for refinement-only balance from balanced input.
Every nonterminal sweep has a birth, so total orchestration costs $O(N_{\rm b})$.

For a gap-independent count, initial and sibling contacts are genealogy
roots with disjoint patch interiors: new sibling patches are internal to
the split cell, whereas inherited patches lie on its boundary.
Replacements partition their parent's patch and strictly increase the
endpoint level sum, bounded by $2dL_{\max}$. Every record has a final
descendant. Disjoint roots give unique ancestry and
$N_{\rm c}^0\le N_{\rm c}^*$, where $N_{\rm c}^*$ counts final contacts.
Final balance gives $N_{\rm c}^*=O_d(n_*)$. Thus
\begin{equation}
\label{eq:contact-final-output-budget}
 N_{\rm c}^0\le N_{\rm c}^*,\qquad
 N_{\rm c}\le(1+2dL_{\max})N_{\rm c}^*,\qquad
 N_{\rm c}^*=O_d(n_*).
\end{equation}
Full face/gluing identities keep these genealogies distinct, including on
self-gluings. Combined with \eqref{eq:incidence-birth-budget}, this gives
$N_{\rm c}=O_d(\min\{2^{g(d-1)},1+L_{\max}\}\,n_*)$.

Each replacement tests at most $4^d$ child pairs. Retired entries are visited
once per endpoint. Sibling creation, mask unions and radix ordering give
\[
\begin{aligned}
 T_{\rm upd}
 &=O_d\!\left((1+w)(N_{\rm b}+N_{\rm c}+N_{\rm t}+n_*)\right)\\
 &=O_d\!\left((1+w)\min\{2^{g(d-1)},1+L_{\max}\}\,n_*\right).
\end{aligned}
\]
Workspace is $O_d(n_1+N_{\rm b}+N_{\rm c})$.
For adapted balanced input, $g\le3$ by \Cref{prop:one-step-imbalance}.
\end{proof}

\subsection{Complexity of adaptive cycles}
\label{sec:composed-accounting}

Let $n,q,n_*$ count input descriptor records, requests and final descriptor
records, and $N_j$ records at selected rank $j\in J$.
Radix sorting and merging resolve flat requests, validate active keys and
apply refinement precedence in $O_{d,w}(n+q)$ work and workspace.

\Needspace{12\baselineskip}
\begin{theorem}[Flat sequential mesh-cycle complexity]
\label{cor:flat-cycle-accounting}
Fix $d,L_{\max}$, a key width $w\ge L_{\max}$ and the macro-face gluings.
For a valid one-level plan on an arbitrary valid flat input forest,
adaptation, balance and optional trace-comparability closure, followed by
materialization of the canonical final forest, input-to-final incidence,
selected rank cuts and incidence between consecutive selected cuts, take work
\begin{equation}
\label{eq:composed-operation-cost}
 T_{\rm cycle}
 =O_{d,w}\!\left(n+q+n_*+\sum_{j\in J}N_j\right).
\end{equation}
\end{theorem}
\begin{proof}
By \Cref{prop:flat-normalization,prop:reduction-work},
\eqref{eq:recursive-contact-work}, \Cref{cor:maintained-incidence} and
\Cref{prop:adaptation-incidence,thm:rank-characterization}, respectively,
\[
 T_{\rm cycle}\le C_{d,w}\left[
 \underbrace{n+q}_{\text{adaptation}}+
 \underbrace{n+n_*}_{\text{closure}}+
 \underbrace{n+n_*+\sum_{j\in J}N_j}_{\text{incidence, rank cuts}}
 \right].
\]
\end{proof}
Reading input and writing complete output costs
$\Omega(n+q+n_*+\sum_jN_j)$, proving optimality.
More explicitly, let $g_{\rm in}$ be the input's maximal tangential gap,
$g=\min\{L_{\max},g_{\rm in}+2\}$ (\Cref{prop:one-step-imbalance}) and
$B_g=2^{g(d-1)}$, $\widehat B_g=\min\{B_g,1+L_{\max}\}$.
Since $N_{\rm c}^0=O_d(n_*)$ and $D\le dL_{\max}\le dw$, recursive routing
contributes only $O_d(n+(1+w)n_*)$.
Combining this with \eqref{eq:contact-final-output-budget} and binary radix ordering gives
\begin{equation}
\label{eq:cycle-gap-depth-cost}
 T_{\rm cycle}=O_d\!\left((1+L_{\max})^{d+1}n+
 (1+w)\left[q+\widehat B_g(n+n_*)+\sum_{j\in J}N_j\right]\right).
\end{equation}
The constant depends only on $d$; face-balanced input has $g\le3$.
Thus depth dependence is polynomial even for unbalanced input.

\section{Conforming finite-element realization}
\label{sec:finite-element-realization}

The mesh algorithms provide active cells, shared-face patches and nested coarse
meshes. Polynomial trace equations give conforming high-order spaces even when
face partitions cross. We construct transfers between these spaces, quantify
the constraints' effect on approximation and obtain mapped heat stability
from cellwise inverse bounds.

\subsection{Conforming finite-element spaces}
\label{sec:conforming-spaces}
Let $\Omega\subset\R^d$ be a bounded Lipschitz domain with a conforming
macro mesh $(M_m)_{m\in\Mcal}$
\cite{ernGuermondFiniteElements2021}*{Chapter~8}. The macro-cell closures cover
$\overline\Omega$ and meet in full common subfaces, with disjoint interiors.
Each bi-Lipschitz map
$G_m:[0,1]^d\to\overline{M_m}$ matches neighbors through the signed permutations
\eqref{eq:signed-permutation}. Refinement may create hanging interfaces.
An active cell $K$ has map $F_K=G_m\circ A_K$,
$A_K(\widehat x)_a=2^{-\ell_a(K)}(i_a(K)+\widehat x_a)$, with reference
widths $\widehat h_{K,a}=2^{-\ell_a(K)}$. For $p\ge1$, we use tensor-product nodes from
any $p+1$ distinct points in $[0,1]$.
Let $\Gamma_D$ be the union of exterior macro faces on which
zero trace is prescribed. Set $\mathbb Q_p=\operatorname{span}\{x^\nu:0\le\nu_a\le p\}$.
The broken and conforming spaces are
\[
 W_h^p=\{v\in L^2(\Omega):v|_K\circ F_K\in\mathbb Q_p\},
 \qquad V_h^p=\{v\in W_h^p\cap H^1(\Omega):v|_{\Gamma_D}=0\}.
\]
The vector $\widetilde u\in\R^{N_h}$ lists the cell nodal coefficients
independently, before constraints are imposed.
The constraint matrix $B_h$ equates both traces at the physical images of
$(p+1)^{d-1}$ interpolation nodes of each overlay patch $H=F\cap G$ in the common reference-face chart,
and sets trace evaluations to zero on such a grid on each active leaf face
contained in $\Gamma_D$.
If $H=F\subseteq G$, the coarse trace is interpolated to the fine trace grid.
These are nodal assignments when endpoints are included. On crossed contacts,
the equations use evaluations on both sides.

For $Q_1$ on the two-strip face in \Cref{fig:trace-overlay}, each vertical
and horizontal strip carries a bilinear polynomial in the common face chart. Agreement on an open
overlay rectangle makes the corresponding polynomials identical. Every
vertical strip meets both horizontal strips, so all four polynomials must
be restrictions of one bilinear polynomial on the whole face. The overlay
specifies where to test this agreement, it adds no FE unknowns.

\begin{proposition}[Conforming realization by trace equations]
\label{thm:conforming-kernel}
With all overlay-patch and Dirichlet rows included, $\ker B_h$ represents exactly
$V_h^p$. This includes unbalanced meshes and crossed traces. For
trace-comparable partitions the one-sided trace equations suffice.
\end{proposition}
\begin{proof}
Let $v\in W_h^p$ have nodal vector $\widetilde u$.
Bi-Lipschitz maps give $v|_K\in H^1(K)$. Coordinatewise affine trace
transitions, up to signed permutation, give $\mathbb Q_p$-traces in the common chart.
Unisolvence gives
\[
 B_h\widetilde u=0
 \quad\Longleftrightarrow\quad
 \left[
 \tr_H(v|_K)=\tr_H(v|_L)\ \forall H,\quad
 v|_{\Gamma_D}=0
 \right].
\]
For $\psi\in C_c^\infty(\Omega)$, with jumps oriented by the unit normal $n$,
cellwise integration by parts gives
\[
 \sum_K\int_K v\,\partial_a\psi
 =-\sum_K\int_K\partial_a v\,\psi
   +\sum_H\int_H[v]\,n_a\psi.
\]
Overlay patches partition the internal interfaces up to null sets, so zero
jumps give the $L^2$ weak gradient. Conversely, an $H^1$ function has equal
traces there. Together with the boundary equations, this proves
$\ker B_h=V_h^p$ in nodal coordinates.
On comparable patches these are precisely the one-sided equations.
\end{proof}

Cellwise directional degrees $\boldsymbol p_K=(p_{K,a})_a\in\N^d$ are also
admissible: on each overlay patch, the trace equations use the componentwise
larger transported degree
(supplement, \emph{Directional polynomial degrees}).

A \emph{constraint basis matrix} $E_h$ contains the broken nodal vectors of
$n_h=\dim V_h^p$ independent conforming basis functions. Its columns span
$\ker B_h$. Multiplication by $E_h$ fills every cell's coefficients from the
free degrees of freedom (DoFs). $E_h^{\mathsf T}$ collects cell residuals
back into the free variables.
For the broken matrix $\widetilde A_h$ and load $\widetilde b_h$,
$\widetilde u=E_hu$, $u\in\R^{n_h}$, gives the constrained system
\begin{equation}
\label{eq:conforming-operator}
 A_h=E_h^{\mathsf T}\widetilde A_hE_h,\qquad
 b_h=E_h^{\mathsf T}\widetilde b_h,\qquad
 n_h=N_h-\operatorname{rank}(B_h).
\end{equation}
On each MPI rank, the production builder uses one-sided rows on
trace-comparable contacts and two-sided rows on crossed contacts. Joint elimination resolves
dependencies across faces before constructing deal.II affine constraints.
The supplement's \emph{Two-sided production constraints and crossed hierarchies} specifies grouping, pivot selection and tolerances.

\subsection{Finite-element transfer}
\label{sec:fe-transfer}
Prolongation represents the same coarse FE function on the fine mesh.
Cell incidence identifies which coarse polynomial to evaluate on each fine
cell, including when traces cross (\Cref{thm:conforming-kernel}).
For same-macro cells $K=G_m(Q_\kappa)$ and
$L=G_m(Q_\eta)$, injectivity gives $K\cap L=G_m(Q_\kappa\cap Q_\eta)$.
Positive-volume overlaps are mapped dyadic cells \eqref{eq:dyadic-dichotomy},
providing cross-mass integration domains for $L^2$ projection between nonnested
adapted spaces \cite{farrellMaddisonConservative2011}. Constraint elimination,
transfer assembly and projection solves are outside the mesh-work bound.

Let $V_j^p$ be the conforming space on $\Fcal_j$ with common $\Gamma_D$.
All levels inherit the same macro maps, so $V_k^p\subseteq V_j^p$ for $j<k$
by \Cref{thm:rank-characterization}: polynomial restriction preserves
$\mathbb Q_p$, continuity and boundary values. With constraint bases $E_j,E_k$
\eqref{eq:conforming-operator}, let $\widetilde P_{j,k}$ evaluate coarse
polynomials at fine nodes. The unique prolongation from rank $k$ to $j$ satisfies
$E_jP_{j,k}=\widetilde P_{j,k}E_k$.
Both sides list the fine-cell nodal values of the same coarse function:
the left uses fine free coefficients, and the right evaluates coarse polynomials.
Residuals return to the coarse level by $P_{j,k}^{\mathsf T}$.
The same identity covers mesh and degree coarsening whenever
$p_{K,a}\ge p_{L,a}$ for every fine cell $K\subseteq L$ and every axis $a$
in their common macro coordinates.

\subsection{Approximation and heat stability}
\label{sec:fem-approximation}
On hanging or crossed faces, continuity determines some nodal values from
others, possibly distant ones. The reconstruction weights and distances
determine the approximation error.
Interpolation uses Cartesian maps
$F_K(\widehat x)_a=b_{K,a}+h_{K,a}\widehat x_a$
on $Q=\Omega_x\times(0,T)$, with rectangular $\Omega_x\subset\R^{d-1}$,
$d\ge2$, physical time $t=x_d$, and zero traces on the lateral and initial
boundary $\Gamma_D$.
\subsubsection*{Constrained interpolation}
The reconstruction uses $p\ge1$ and fixed reference interpolation nodes
including $0,1$.
Let $z_i$ be the broken cell nodes and $\overline B_h$ the interior trace
rows of $B_h$ (\Cref{thm:conforming-kernel}). We select nodal samples that
uniquely determine a function in $\ker\overline B_h$: a maximal independent
set of evaluations on $\Gamma_D$, extended by other independent evaluations.
Let $\mathcal Z$ index these samples.
The nodal reconstruction matrix $\mathsf V\in\R^{N_h\times|\mathcal Z|}$
recovers all nodal values:
\begin{equation}
\label{eq:constrained-interpolant}
 \overline B_h\mathsf V=0,\qquad
 \mathsf V_{\mathcal Z,:}=I,\qquad
 (I_hu)(z_i)=\sum_{j\in\mathcal Z}\mathsf V_{ij}u(z_j).
\end{equation}
Only interior continuity is imposed in constructing $\mathsf V$, so $I_h$
reproduces every global $\mathbb Q_p$ polynomial. Boundary-first sampling
preserves zero data on $\Gamma_D$ because the selected evaluations span
all boundary evaluations. Terminal time remains free.
For cell $K$ with widths $h_{K,a}$ and node indices $\mathcal N_K$, the
\emph{determining patch} $\omega_K$ is the smallest rectangle containing $K$
and all samples $z_j$ with $\mathsf V_{ij}\ne0$ for some $i\in\mathcal N_K$.
Its widths and maximal row sum measure reach and amplification:
\[
 R_{K,a}=\operatorname{width}_a(\omega_K),\qquad
 \Lambda_K=\max_{i\in\mathcal N_K}\sum_j|\mathsf V_{ij}|,\qquad
  h_K^{-\beta}=\prod_a h_{K,a}^{-\beta_a}.
\]
For $u\in C(\overline Q)$, the continuity constraints give
\begin{equation}
\label{eq:local-best-approximation}
 \|u-I_hu\|_{L^\infty(K)}
 \le C_{d,p}(1+\Lambda_K)
       \inf_{q\in\mathbb Q_p}\|u-q\|_{L^\infty(\omega_K)}.
\end{equation}
Indeed, \eqref{eq:constrained-interpolant} gives
$\|I_h(u-q)\|_{L^\infty(K)}\le
C_{d,p}\Lambda_K\|u-q\|_{L^\infty(\omega_K)}$
for $q\in\mathbb Q_p$, and ${I_hq=q}$.
To account for weak influence from distant samples, define the moments
\[
 M_i^\alpha=\sum_{j\in\mathcal Z}|\mathsf V_{ij}|\,
             |(z_j-z_i)^\alpha|\qquad(i\in\mathcal N_K).
\]

\Needspace{12\baselineskip}
\begin{proposition}[Moment-dependent directional approximation]
\label{prop:constrained-approximation}
For smooth $u$ and $\beta\in\{0,\ldots,p+1\}^d$, with coordinate unit vectors $e_a$,
including on meshes with crossed traces,
\begin{equation}
\label{eq:constrained-approximation}
\begin{split}
 \|\partial^\beta(u-I_hu)\|_{L^\infty(K)}
 \le C_{d,p}\biggl[&\sum_a h_{K,a}^{p+1-\beta_a}
 \|\partial_a^{p+1}\partial^{\beta-\beta_ae_a}u\|_{L^\infty(K)}\\
 &+h_K^{-\beta}\max_{i\in\mathcal N_K}\sum_{|\alpha|=p+1}
 \frac{M_i^\alpha}{\alpha!}
 \|\partial^\alpha u\|_{L^\infty(\omega_K)}\biggr].
\end{split}
\end{equation}
The constant $C_{d,p}$ is independent of aspect ratios and coordinate levels.
\end{proposition}
\begin{proof}
Let $T_i^pu$ be the degree-$p$ Taylor polynomial at $z_i$.
Polynomial reproduction and Taylor's remainder on segments in $\omega_K$ give
\[
 |(I_hu-u)(z_i)|
 =\Big|\sum_j\mathsf V_{ij}(u-T_i^pu)(z_j)\Big|
 \le\sum_{|\alpha|=p+1}\frac{M_i^\alpha}{\alpha!}
       \|\partial^\alpha u\|_{L^\infty(\omega_K)}.
\]
For the cellwise tensor interpolant $\mathcal J_K$, standard interpolation
bounds $\partial^\beta(u-\mathcal J_Ku)$ by the first sum
\cite{apelAnisotropicFiniteElements1999}. Scaled nodal stability bounds
$\partial^\beta(\mathcal J_Ku-I_hu)$ by the second.
\end{proof}

For nonnegative $Q_1$ rows, $M_i^{2e_a}$ is the weighted variance in axis $a$.
Reproduction of $x_ax_b$, $a\ne b$, cancels mixed covariances, so every
quadratic polynomial $q$ satisfies
\[
 (I_hq-q)(z_i)=\tfrac12\sum_a\partial_a^2q\,M_i^{2e_a}.
\]

\begin{proposition}[Logarithmic interpolation loss on balanced meshes]
\label{prop:local-interpolation-loss}
For each integer $2\le r<L_{\max}$, one-level refinement-and-balance cycles from the uniform
four-cell mesh of $Q=(0,1)^2$ produce a face-balanced, trace-comparable mesh
and a square $K_r$ of side $h_r=2^{-r-1}$. For recursive $Q_1$ coarse-edge
interpolation, with unconstrained vertices as samples in
\eqref{eq:constrained-interpolant}, the function $u(x,t)=x(1-x)t^2$ satisfies
\begin{equation}
\label{eq:local-interpolation-loss}
 \frac r8h_r^2\le\|u-I_hu\|_{L^\infty(K_r)}\le rh_r^2,
 \qquad \Lambda_{K_r}=1.
\end{equation}
\end{proposition}
The supplement's \emph{Repeated corner refinement} gives the proof.
At the lower-right vertex $z_i$, $M_i^{2e_1}=rh_r^2/2$ and $M_i^{2e_2}=rh_r^2$:
the normalized moments $M_i^{2e_a}/h_r^2$ grow like $r$, while the patch ratios
$R_{K_r,1}/h_r=2^{r-1}$ and $R_{K_r,2}/h_r=2^r$ grow exponentially.
Face balance does not bound the patch ratio assumed in
\cite{stevensonStornInterpolation2023}*{Sec.~6}.

\subsubsection*{Mapped residual heat equation}
We now return to the mapped conforming meshes of \Cref{sec:conforming-spaces}
and use the residual form of \cite{langerSchafelnerHexahedral2022},
including tilted and crossed interfaces. Maps and inverses have bounded second
derivatives. Here $Q=\Omega_x\times(0,T)$, $\Omega_x\subset\R^{d_x}$, or
$Q=\Omega_x\times(0,T)\times\Upsilon$ with passive parameters
$\mu\in\Upsilon\subset\R^{d-d_x-1}$. Derivatives use physical coordinates $(x,t,\mu)$;
$d=d_x+1$ without parameters.
The space $V_h^p$ has zero initial and spatial lateral traces.
The terminal time and parameter faces are free.
Let $k=k(x,t,\mu)\in\R^{d_x\times d_x}$ be symmetric, uniformly positive and bounded, with cellwise
bounded spatial derivatives. Write
$\mathcal L_k w=\operatorname{div}_x(k\nabla_xw)$ and
$\|\nabla_xw\|_{k,D}^2=(k\nabla_xw,\nabla_xw)_D$.
We consider $u_t-\mathcal L_ku=f$, $f\in L^2(Q)$, with these homogeneous data.
For cellwise $\tau_K>0$, set
\begin{equation}
\label{eq:mapped-heat-residual}
\begin{aligned}
 b_\tau(w,v)&=(w_t,v)_Q+(k\nabla_xw,\nabla_xv)_Q
       +\sum_K\tau_K(w_t-\mathcal L_kw,v_t)_K,\\
 \ell_\tau(v)&=(f,v)_Q+\sum_K\tau_K(f,v_t)_K.
\end{aligned}
\end{equation}
The discrete solution $u_h\in V_h^p$ satisfies
$b_\tau(u_h,v_h)=\ell_\tau(v_h)$ for every $v_h\in V_h^p$.
In $b_\tau(v_h,v_h)$, stabilization contributes
$\tau_K\|v_{h,t}\|_K^2-\tau_K(\mathcal L_kv_h,v_{h,t})_K$ on each cell.
A local inverse bound controls the competing cross term.
Let $C_K$ be the least constant in
$\|\mathcal L_kv_h\|_K\le C_K\|\nabla_xv_h\|_{k,K}$ for mapped local
$\mathbb Q_p$ polynomials. It is finite by finite dimensionality:
$\nabla_xv_h=0$ implies $\mathcal L_kv_h=0$. The supplement gives its
cell-matrix formula, retaining curvature, anisotropy and coefficient variation.

\Needspace{8\baselineskip}
The energy norm measures spatial diffusion, time variation and the terminal
trace. Its companion controls the form after integration by parts in time.
For $w,w_t,\nabla_xw\in L^2(Q)$ with cellwise $\mathcal L_kw\in L^2(K)$, define
\begin{equation}
\label{eq:heat-energy-norms}
\begin{aligned}
 \|w\|_h^2&=\|\nabla_xw\|_{k,Q}^2+
       \sum_K\tau_K\|w_t\|_K^2+\tfrac12\|w(T)\|^2,\\
 \|w\|_{h,*}^2&=\|\nabla_xw\|_{k,Q}^2+
       \sum_K\tau_K\|w_t-\mathcal L_kw-\tau_K^{-1}w\|_K^2
       +2\|w(T)\|^2.
\end{aligned}
\end{equation}

\Needspace{12\baselineskip}
\begin{proposition}[Mapped residual stability and approximation]
\label{prop:heat-energy}
Suppose $\theta<2$, with $\theta=\max_K\sqrt{\tau_K}C_K$, and set $\alpha=1-\theta/2$.
With exact integration, \eqref{eq:mapped-heat-residual} has a unique discrete
solution and, for the homogeneous-data PDE solution $u$ with
$u,u_t,\nabla_xu\in L^2(Q)$ and cellwise $\mathcal L_ku\in L^2(K)$,
\begin{equation}
\label{eq:heat-best-approximation}
\begin{aligned}
 b_\tau(v_h,v_h)&\ge\alpha\|v_h\|_h^2,\\
 \|u-u_h\|_h&\le\inf_{v_h\in V_h^p}
       \bigl(\|u-v_h\|_h+\alpha^{-1}\|u-v_h\|_{h,*}\bigr).
\end{aligned}
\end{equation}
In particular, $\tau_KC_K^2\le1$ gives $\alpha\ge1/2$.
\end{proposition}
\begin{proof}
See the supplement's \emph{Mapped residual stability}.
\end{proof}

For $Q=Q_{xt}\times\Upsilon$, $Q_{xt}=\Omega_x\times(0,T)$, the first-derivative
assumption is $u\in L^2(\Upsilon;H^1(Q_{xt}))$. Terminal traces lie in
$L^2(\Omega_x\times\Upsilon)$. Initial and lateral conditions hold for almost
every parameter. All products integrate over $\Upsilon$.
Square-integrable parameter dependence, including jumps, is allowed:
no parameter derivatives enter the argument
\cite{daworStochasticGalerkin2026}*{Section~2.2.2}.

With positive operator quadrature, coercivity uses quadrature inverse bounds
and the identity $(v_{h,t},v_h)_{Q,q}=\tfrac12\|v_h(T)\|^2$ for $v_h\in V_h^p$.
For multilinear maps, its volume and boundary pullbacks have coordinate
degree at most $2p+d-2$. Thus $q\ge p+\lceil(d-1)/2\rceil$ Gauss points
preserve this identity exactly. The reported $q=p+2$ rule suffices for every
degree in 4D/5D. The supplement gives the degree count and inverse diagnostics.
Diffusion/residual quadrature consistency remains separate.

On Cartesian meshes, \eqref{eq:constrained-approximation} controls the residual
norm's second spatial derivatives. For scalar diffusion, the supplement bounds
$C_K$ using spatial widths and coefficient bounds. Its \emph{One-generation
reconstruction} bounds $\Lambda_K,R_{K,a}/h_{K,a}$ for a common split mask
on a face-to-face mesh, with $\Lambda_K=1$ for $Q_1$.

\subsection{Matrix-free multigrid solves}
\label{sec:multigrid-solver}
The implementation supplies cells, DoFs and geometry to deal.II, patched
for high dimensions \cites{arndtDealiiDesign2021,bangerthLibrary2007}.
Its matrix-free kernels apply cell operators without storing the global matrix,
using sum factorization, vectorized batching and distributed constraint handling
\cite{kronbichlerKormannMatrixfree2012}.

Geometric multigrid uses rank cuts \eqref{eq:active-rank-set} at fixed degree.
Hybrid multigrid \cite{mitchellHpMultigrid2010} first coarsens the degree by
$p_{r+1}=\max(2,\lfloor p_r/2\rfloor)$ to $p_r\le2$, then uses rank cuts.
Both use FE
prolongation and transpose restriction (\Cref{sec:fe-transfer}) with deal.II's
V-cycle, Chebyshev and Krylov components \cite{munchMultigrid2023}.
Each level rediscretizes the same form. Solver settings are in
\Cref{tab:multigrid-solvers,tab:heat-solvers} and the supplement's
\emph{Measured solver settings}.

The heat smoother uses short residual solves for the full level operator,
including its nonsymmetric terms. These solves are preconditioned by
Chebyshev on the shifted energy operator
\begin{equation}
\label{eq:heat-smoother}
 H_\ell=\nabla_x^*k\nabla_x+\tau\partial_t^*\partial_t+\alpha_\ell I,
 \qquad \alpha_\ell>0.
\end{equation}
Stars denote $L^2$ adjoints.
In the conforming basis $\{\varphi_{\ell,i}\}$, the discrete shift is
$\alpha_\ell M_\ell$, where $(M_\ell)_{ij}=(\varphi_{\ell,j},\varphi_{\ell,i})$
is evaluated with the operator quadrature.
In the moving-layer examples, $\alpha_\ell=1$ and every level retains the
finest $\tau$ and full stabilized form. Flexible GMRES (FGMRES) \cite{saadFlexible1993}
accommodates the residual-dependent inner solves.
Only the small coarse system is assembled.

\section{Numerical experiments}
\label{sec:numerical-experiments}

The experiments test the complete construction in four and five dimensions:
mesh cycles at millions of cells, directional accuracy gains with
$Q_2$--$Q_4$ elements, and coarse corrections from generated hierarchies.
The accuracy plots show all 524 samples at their attained numbers of free
degrees of freedom (DoFs); every error curve decreases.
FE comparisons use two MPI ranks on distinct cores. Each solver pair uses
the same discrete system. Timings are three-run medians.
Multigrid setup includes rank cuts, coarse FE spaces, operators, transfers and smoothers
(supplement, \emph{Measured solver settings}).
The supplement and accompanying reproduction package give implementation checks,
sources, deal.II patches, measurements and reproduction commands for the
reported experiments.

\subsection{Cost of mesh cycles}
\label{sec:algorithm-scaling}
\begin{figure}[!htbp]
\centering
\includegraphics[width=\linewidth]{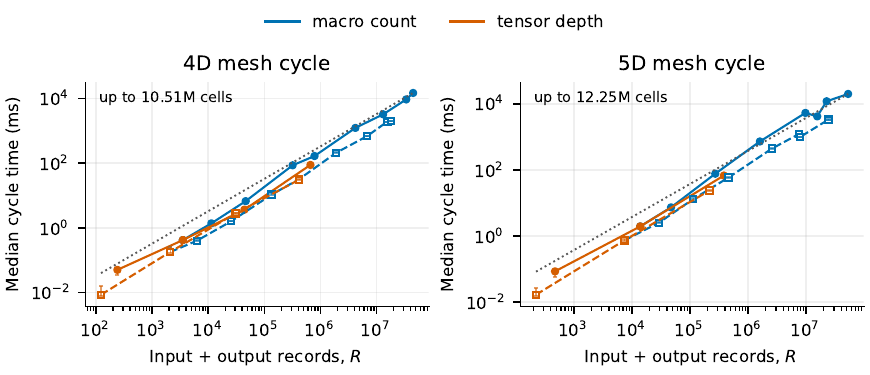}
\caption{Complete flat mesh cycles in 4D (left) and 5D (right), including
selected hierarchies and incidence. Curves show five-run medians/ranges
(three for the largest cases). Solid/dashed curves refine/coarsen.
Dotted lines have slope one.}
\label{fig:mesh-cycle-cost}
\end{figure}

The \emph{accounting size} $R=n+q+n_*+\sum_jN_j$ counts input, request,
final and selected hierarchy records (\Cref{cor:flat-cycle-accounting}).
Flat cycles reach $10.51$ million leaves in 4D and $12.25$ million in 5D
(\Cref{fig:mesh-cycle-cost,tab:large-cycles}). Their one-core median times are $14.76$ and
$19.72$ seconds, including adaptation, both closures, topology, rank cuts
and incidence. The respective accounting sizes are $45.02$ and $52.64$ million.

Macro-growth fixtures use level-two tensor seeds. Tensor-depth fixtures vary
the uniform seed depth within one macro. The supplement includes the corner-depth family
and implementation and timing details.

\begin{table}[!htbp]
\centering
\small
\caption{Complete cycles with millions of cells: three-run medians.
Input records $n$, final leaves $L_*$ and accounting size $R$ are in millions.
The closure factor $\gamma=n_*/n_{\rm ad}$ compares descriptor counts after
and before both closures.}
\label{tab:large-cycles}
\begin{tabular}{llrrrrrr}
\hline
$d$ & plan & $n/10^6$ & $L_*/10^6$ & $\gamma$ & $R/10^6$ & seconds & $\mu$s$/R$ \\
\hline
4 & refine & 4.885 & 10.511 & 1.95 & 45.025 & 14.76 & 0.328 \\
4 & coarsen & 4.885 & 2.248 & 1.04 & 18.192 & 2.05 & 0.113 \\
5 & refine & 6.780 & 12.247 & 1.61 & 52.645 & 19.72 & 0.375 \\
5 & coarsen & 6.780 & 2.952 & 1.03 & 24.149 & 3.30 & 0.137 \\
\hline
\end{tabular}
\end{table}

Across all fixtures, median cycle time per accounting record is
$0.040$--$0.328\,\mu\mathrm{s}$ in 4D and $0.049$--$0.544\,\mu\mathrm{s}$ in 5D.
The 5D 324-macro fixture has $\gamma=3.28$ and 19 cumulative balance sweeps
across trace-comparability iterations. The 1280-macro fixture has $\gamma=1.39$
and three.
Medians follow warm-up on a pinned Intel Core Ultra~7 155U core.

\subsection{Mapped Poisson}
\label{sec:fem-verification}
\label{sec:solver-comparison}

\begin{figure}[!htbp]
\centering
\includegraphics[width=\linewidth]{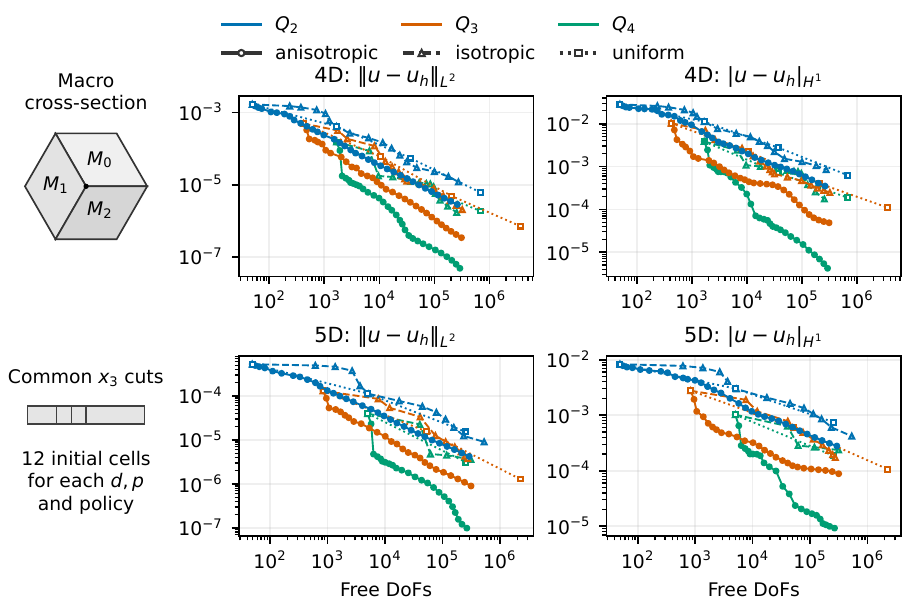}
\caption{Directional accuracy gains for mapped Poisson in 4D (top) and 5D
(bottom). Left: common geometry and initial mesh. Center/right: $L^2$ and
physical $H^1$ seminorm errors. Both adaptive policies use exact-solution scores.}
\label{fig:poisson-refinement}
\end{figure}

We compare split directions for a narrow layer on mapped meshes with cyclic
macro connectivity (\Cref{fig:poisson-refinement}).
We solve $-\Delta u=f$ with zero Dirichlet data on $H_{\rm hex}\times(0,1)^{d-2}$, $d=4,5$,
where $H_{\rm hex}$ is the regular hexagon with vertices
$v_j=(\cos(j\pi/3),\sin(j\pi/3))$, $j=0,\ldots,5$.
Three affine macros form a cycle around the central junction:
\[
 G_m(\xi)=(\xi_1v_{2m}+\xi_2v_{2m+2},\xi_3,\ldots,\xi_d),\qquad m=0,1,2
 \quad\text{(indices modulo six)}.
\]
Their Jacobians are non-diagonal. The exact solution is
\[
 u=e^{c\cdot x-((x_3-.37)/.12)^2}\prod_{a=3}^d x_a(1-x_a)
       \prod_{j=0}^2\bigl[3/4-(n_j\cdot(x_1,x_2))^2\bigr],
\]
with $c$ the first $d$ entries of $(.2,.35,.4,.6,.7)$ and
$n_j=(\cos(\pi/6+j\pi/3),\sin(\pi/6+j\pi/3))$.
The polynomial factors impose zero boundary data. The Gaussian concentrates
variation near $x_3=.37$. Analytic derivatives give $f=-\Delta u$.
All runs share twelve initial cells, with $x_3$ cuts
$0,1/4,3/8,1/2,1$ tracking this layer.

We use the conforming spaces of \Cref{thm:conforming-kernel}, with
Gauss-Lobatto interpolation nodes, $p+1$ Gauss points per axis for operators
and ten for loads/errors.
The refinement curves use multigrid-preconditioned FGMRES at tolerance $10^{-11}$.

\subsubsection*{Refinement and accuracy.}
Exact-solution scores isolate directional splitting from error-estimator quality.
Directional interpolation motivates heuristic Taylor scores for Poisson and
heat, using reference cells $Q_\kappa$ and widths $\widehat h_{K,a}$
(\Cref{sec:conforming-spaces}):
\begin{equation}
\label{eq:taylor-score}
 s_{K,a}=\int_{Q_\kappa}\left[
 \frac{\widehat h_{K,a}^{p+1}}{(p+1)!}\partial_{\xi_a}^{p+1}(u\circ G_m)(\xi)
 \right]^2|\det DG_m(\xi)|\,d\xi,
 \quad S_K=\sum_a s_{K,a}.
\end{equation}
Three-point tensor Gauss quadrature evaluates them. Both policies mark cells
accounting for at least 50\% of the total score $\sum_K S_K$, with key tie-breaking.
Anisotropic requests (A) split the largest-score axis, ties to the smallest
axis; isotropic requests (I) split all $d$ axes.
Balance closure and trace-comparability closure
(\Cref{thm:balance-closure,thm:trace-comparability-closure}) follow each step,
permitting one-sided constraints on the finest mesh.
Uniform controls bisect every cell in all axes.
Each refinement sequence stops at its first sample with at least $250{,}000$ free DoFs.

At $50{,}000$ free DoFs, the isotropic $L^2$ error is $2.8$--$42.3$ times
the anisotropic error in 4D and $4.6$--$9.4$ times in 5D across $Q_2$--$Q_4$.
These values use log-log interpolation between bracketing samples.
Final anisotropic $Q_4$ samples attain the smallest errors
(\Cref{fig:poisson-refinement}). Uniform controls reach $3.62$ million free DoFs.
For 4D $Q_4$, higher score quadrature gives an error ratio of $42.73$ (supplement).
The coarse uniform $Q_3/Q_4$ samples are preasymptotic on the layer.

\subsubsection*{Multigrid performance.}
The generated rank hierarchies reduce outer FGMRES iterations from
141--821 to 7--29
(\Cref{tab:multigrid-solvers}). Median setup-plus-solve time is lower in
six of eight cases, by up to a factor of $4.54$.
For $Q_2$ and $Q_4$, each pair uses the first anisotropic or isotropic state
at or above $100{,}000$ free DoFs in its dimension. Both solvers use FGMRES,
with single-level Jacobi ($B$) or geometric multigrid ($h$,
\Cref{sec:multigrid-solver}).
Every elementary rank is used, including crossed coarse traces with the
two-sided constraints and transfer of \Cref{sec:fe-transfer}.

\begin{table}[!htbp]
\centering
\small
\caption{Poisson: generated hierarchies reduce iterations in every pair.
Jacobi ($B$), geometric multigrid ($h$); A/I: anisotropic/isotropic.
Free DoFs, outer FGMRES iterations and three-run median setup-plus-solve
times on two MPI ranks (\Cref{sec:numerical-experiments}).}
\label{tab:multigrid-solvers}
\begin{tabular}{lllrrrrrrr}
\hline
\multirow{2}{*}{$d$} & \multirow{2}{*}{Degree} & \multirow{2}{*}{Mark} & \multirow{2}{*}{DoFs} & \multirow{2}{*}{Levels} & \multicolumn{2}{c}{Iterations} & \multicolumn{2}{c}{Time (s)} & \multirow{2}{*}{$t_B/t_{h}$}\\
\cline{6-7}\cline{8-9}
 & & & & & $B$ & $h$ & $B$ & $h$ & \\
\hline
\multirow{4}{*}{4} & \multirow{2}{*}{$Q_2$} & A & 115{,}029 & 9 & 550 & 23 & 7.485 & 4.513 & 1.66 \\
 &  & I & 180{,}337 & 7 & 475 & 15 & 8.942 & 2.897 & 3.09 \\
\cline{2-10}
 & \multirow{2}{*}{$Q_4$} & A & 104{,}245 & 7 & 821 & 29 & 6.402 & 2.580 & 2.48 \\
 &  & I & 104{,}097 & 5 & 272 & 7 & 1.560 & 0.344 & 4.54 \\
\hline
\multirow{4}{*}{5} & \multirow{2}{*}{$Q_2$} & A & 120{,}433 & 9 & 253 & 19 & 8.681 & 11.814 & 0.73 \\
 &  & I & 113{,}443 & 6 & 212 & 17 & 3.893 & 4.013 & 0.97 \\
\cline{2-10}
 & \multirow{2}{*}{$Q_4$} & A & 107{,}747 & 8 & 329 & 18 & 4.807 & 3.272 & 1.47 \\
 &  & I & 118{,}931 & 4 & 141 & 9 & 1.716 & 0.756 & 2.27 \\
\hline
\end{tabular}
\end{table}

\subsection{Moving heat layers}
\label{sec:marking-comparison}
\begin{figure}[!htbp]
\centering
\includegraphics[width=\linewidth]{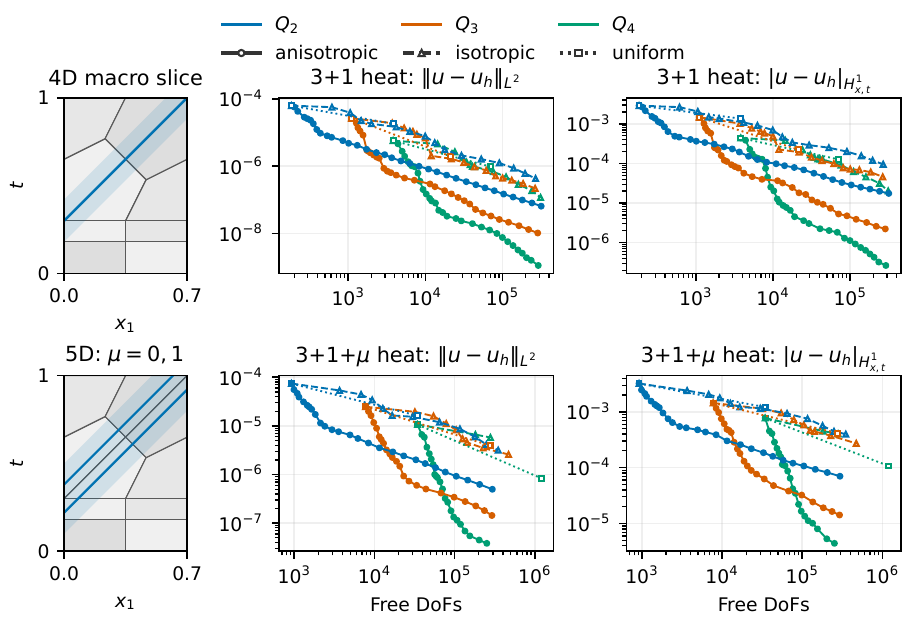}
\caption{Directional refinement resolves moving space-time (top) and
parametric heat layers (bottom) with fewer free DoFs.
Left: macro-mesh sections in $(x_1,t)$ at $x_2=.12$, $x_3=.18$
(bottom: $\mu=0,1$). Gray lines show macro boundaries. Blue lines mark
layer centers ($\phi=0$).
Center/right: $L^2$ errors and space-time $H^1$ seminorm errors.
Adaptive curves use exact-solution Taylor scores.}
\label{fig:marking-comparison}
\end{figure}

The heat tests adapt jointly in space and time to resolve a moving layer
without time slabs (\Cref{fig:marking-comparison}). The parametric extension also resolves changes in
diffusivity and layer position on a single joint mesh.
Let $\Omega_x=\{x\in\mathbb R^3:x_i>0,\ \sum_i x_i<1\}$,
$\mathcal Q_4=\Omega_x\times(0,1)_t$, and
$\mathcal Q_5=\mathcal Q_4\times(0,1)_\mu$.
On both domains, we solve $u_t-\nabla_x\cdot(k\nabla_xu)=f$, with
\[
 u=t(1+.2t)x_1x_2x_3(1-x_1-x_2-x_3)\,
 e^{c\cdot z-(\phi/.12)^2},\qquad c_j=.10+.03(j-1).
\]
In 4D, $z=(x,t)$, $k=1$, and $\phi=t-\sum_i x_i$.
In 5D, $z=(x,t,\mu)$ and
\[
 k=.2+1.8\mu^2+.3\mu(x_1+tx_2),\qquad
 \phi=t-\sum_i x_i-.16(\mu-.5),\qquad .2\le k\le2.3.
\]
The load is $f=u_t-\nabla_x\cdot(k\nabla_xu)$.
The Gaussian factor is centered on the moving surface $\phi=0$. Varying the parameter
shifts this surface and changes diffusivity. Lateral and initial traces vanish,
while terminal time and both parameter boundaries are free. $\mu$ enters as
a discretized coordinate with no derivatives in the PDE.

Simplex subdivisions give 20 mapped hypercube macros on $\mathcal Q_4$. Two parameter
bands with shared control vertices warped in time and $\mu$ give 40 multilinearly
mapped macros on the fixed domain $\mathcal Q_5$. Their geometry is neither a time nor
a parameter extrusion, and the central layer at $\mu=.5$ follows macro faces.
All degrees and refinement policies share these coarse meshes, whose warp
and bi-Lipschitz certificate are given in the supplement for the map
assumptions of \Cref{thm:conforming-kernel}.

The experiments use \eqref{eq:mapped-heat-residual} with a common
$\tau_K=\tau=.1\rho_{\min}^2/[k_{\max}(p+1)^4]$ and
 $\rho_{\min}^{-1}=\max_{K,v}\|J_K(v)^{-1}P_{xt}\|_F$, with $v$ ranging over reference vertices.
Here $k_{\max}=1$ or $2.3$,
$J_K(\widehat x)=DG_m(A_K(\widehat x))\operatorname{diag}(\widehat h_K)$, and
$P_{xt}:\R^4\to\R^d$ is the identity for $d=4$ and appends a zero parameter
component for $d=5$. The vertex scale sets $\tau$. Stability requires the
inverse bound of \Cref{prop:heat-energy}.
With exact integration, a global $\tau$ satisfying $\theta<2$ on the finest
level satisfies it on every coarser level, since local polynomial inclusion gives
$C_P\le\max_{K\subseteq P}C_K$ for each coarse cell $P$
(supplement, \emph{Mapped residual stability}).
Quadrature diagnostics on the initial $Q_2$ macros give
numerical maxima $\tau C_{K,q}^2<0.0394$ in 4D and $<0.0338$ in 5D
(supplement). These checks do not certify stability along
the full adaptive trajectories.
\subsubsection*{Refinement and accuracy.}
Local refinement in \Cref{fig:marking-comparison} uses the $Q_2$--$Q_4$ scores
\eqref{eq:taylor-score} and the quadrature, marking, closures and
$250{,}000$-DoF stopping rule of \Cref{sec:fem-verification}.
Uniform controls refine the initial mesh globally once.
On these layer-fitted macros, a one-axis split adds resolution across the layer
with two children. An all-axis split creates $2^d$.

We solve \eqref{eq:mapped-heat-residual} with the hybrid multigrid solver of
\Cref{sec:multigrid-solver} at FGMRES tolerance $10^{-10}$.
Operator quadrature uses $p+2$ Gauss points per
axis and loads/errors use eight, with curvature terms included in physical
second derivatives. Errors integrate over all coordinates, while the
space-time $H^1$ seminorm uses $(\nabla_x,\partial_t)$, excluding parameter derivatives.

At $50{,}000$ free DoFs, isotropic $L^2$ errors are $6.4$--$54.3$ times the
anisotropic errors for space-time heat and $4.2$--$34.5$ times for parametric
heat across $Q_2$--$Q_4$. These values use log-log interpolation between
bracketing samples.
Final anisotropic $Q_4$ samples attain the smallest errors in both problems,
although higher degree does not improve every intermediate DoF budget.
The score-quadrature and stabilization controls also give ratios above one
(supplement).

\subsubsection*{Multigrid performance.}
For the nonsymmetric heat operator, hybrid multigrid reduces outer iterations
from 25--140 to 6--37
(\Cref{tab:heat-solvers}). Median setup-plus-solve time is lower in one of
eight cases. Each residual-smoothing application contains several inner
operator applications, so fewer outer iterations need not reduce time.
Solve time exceeds setup time in all eight heat cases. The supplement's
\emph{Iteration histories and timing spread} gives the full refinement-curve
iteration histories, separate setup/solve times and paired three-run ratios.

The $Q_2/Q_4$ samples are selected as for Poisson. The baseline applies eight
Chebyshev steps to $H_0$ (\eqref{eq:heat-smoother}, $\alpha_0=1$).
The hierarchy uses leaf-halving cuts, including crossed traces
(\Cref{sec:hierarchy}); ``Levels'' counts geometric and polynomial levels.
The supplement lists the ranks.

\begin{table}[!htbp]
\centering
\small
\caption{Heat: generated hierarchies reduce iterations in every pair.
Baseline ($B$), hybrid multigrid ($hp$); A/I: anisotropic/isotropic.
Free DoFs, outer FGMRES iterations and three-run median setup-plus-solve
times on two MPI ranks (\Cref{sec:numerical-experiments}).}
\label{tab:heat-solvers}
\begin{tabular}{lllrrrrrrr}
\hline
\multirow{2}{*}{Problem} & \multirow{2}{*}{Degree} & \multirow{2}{*}{Mark} & \multirow{2}{*}{DoFs} & \multirow{2}{*}{Levels} & \multicolumn{2}{c}{Iterations} & \multicolumn{2}{c}{Time (s)} & \multirow{2}{*}{$t_B/t_{hp}$}\\
\cline{6-7}\cline{8-9}
 & & & & & $B$ & $hp$ & $B$ & $hp$ & \\
\hline
\multirow{4}{*}{$3+1$} & \multirow{2}{*}{$Q_2$} & A & 103{,}222 & 6 & 82 & 12 & 21.425 & 18.281 & 1.17 \\
 &  & I & 139{,}320 & 4 & 37 & 8 & 10.791 & 11.120 & 0.97 \\
\cline{2-10}
 & \multirow{2}{*}{$Q_4$} & A & 116{,}468 & 6 & 111 & 22 & 10.723 & 10.792 & 0.99 \\
 &  & I & 110{,}944 & 4 & 45 & 12 & 4.318 & 5.630 & 0.77 \\
\hline
\multirow{4}{*}{$3+1+\mu$} & \multirow{2}{*}{$Q_2$} & A & 104{,}622 & 6 & 71 & 16 & 61.898 & 90.401 & 0.68 \\
 &  & I & 116{,}452 & 3 & 25 & 6 & 22.758 & 27.702 & 0.82 \\
\cline{2-10}
 & \multirow{2}{*}{$Q_4$} & A & 107{,}260 & 4 & 140 & 37 & 26.088 & 36.858 & 0.71 \\
 &  & I & 181{,}532 & 3 & 56 & 15 & 23.085 & 29.513 & 0.78 \\
\hline
\end{tabular}
\end{table}

\FloatBarrier
\section{Conclusion}
\label{sec:conclusion}

Canonical representation, output-optimal cycles and trace equations connect
directional meshes to conforming finite elements. At fixed $d,w$, cycle work
has polynomial depth dependence. From balanced input, refinement-only balance
stays on sheets, with an attained geometric bound. Reconstruction moments
bound interpolation errors and capture the logarithmic loss that nodal
interpolation can suffer on balanced meshes.
On layer-fitted 4D/5D macro meshes, directional refinement gains accuracy and
generated hierarchies reduce outer iterations in all 16 tested solver
configurations and setup-plus-solve time in seven.

Uniform local quasi-interpolation on balanced meshes with nested traces and
depth-independent completion by normal grading remain open. Both arise in
anisotropic dual-weighted residual adaptivity
\cites{bauseAnisotropicDwr2026,margenbergAnisotropicHp2026} for
space-time optimal control \cite{endt} and stochastic Galerkin
methods \cite{daworStochasticGalerkin2026}.

\FloatBarrier
\appendix
\section{Auxiliary proofs and implementation bounds}
\label{sec:implementation}

\begin{lemma}[Work of a section or fusion]
\label{lem:section-work}
Section and fusion \eqref{eq:canonical-section}--\eqref{eq:canonical-fusion} take
$O_d(\sum_{\eta\in\Lcal}[1+\|\ell(\eta)-\ell(\kappa)\|_1])$
work and allocated records on canonical inputs with participating leaves
$\Lcal$ in region $\kappa$, reusing unchanged descriptors by handle.
\end{lemma}
\begin{proof}
Pair each section invoking fusion with that fusion; both return the same region.
Child queries have proper regions with disjoint leaf sets. Thus each leaf
$\eta$ is charged along a strict chain
$\kappa\preceq\kappa_0\prec\cdots\prec\kappa_q\preceq\eta$, with
$q+1\le1+\|\ell(\eta)-\ell(\kappa)\|_1$.
Each region needs $O_d(1)$ work and records; direct selections return handles.
Summing these charges proves both bounds.
\end{proof}

\subsection{Branching of refinement propagation}
\label{app:propagation-branching}
Fix $(K,a)$ in the graph of \Cref{thm:marked-growth} and a face $f$ normal
to $n\ne a$, using transported axes \eqref{eq:comparison-axis-maps}.
An outgoing neighbor $L$ contains $K$'s $a$-interval. Balance and dyadic
nesting make its projection onto the other $d-2$ coordinates a nonempty
union of a subset of the $q=2^{d-2}$ half-grid boxes. Disjoint projections
give at most $q$ contacts per face, hence $\deg^+(K,a)\le2(d-1)q=D_d$.
For an incoming neighbor $P$, choose a box in its projection. That box and
its $d-2$ face-adjacent boxes are excluded: boxes inside $P$'s projection
by disjoint interiors. An adjacent box outside it would make $P$ and $L$
share a face, since their normal intervals start at $f$ and the remaining
intervals overlap. This contradicts balance:
$\ell_a(P)-\ell_a(L)=(\ell_a(K)+1)-(\ell_a(K)-1)=2$.
The arrival face loses at least $d-1$ boxes, giving
\[
 \deg^+(K,a)\le2(d-1)q-(d-1)=B_d.
\]
The count includes parallel contacts and self-gluings.

\paragraph*{\textbf{Acknowledgement}}
\begingroup\emergencystretch=1em
This work was supported by FEM POWER II (grant ZS2024/\allowbreak06/\allowbreak18815), funded by the
European Social Fund Plus (ESF+).
\par\endgroup

\paragraph*{\textbf{Code and data availability}}
Code, data and reproduction materials: \cite{margenbergOmnitreeArchive2026}.

\paragraph*{\textbf{Declaration}}
The author used \mbox{ChatGPT} for language editing and consistency
checks, reviewed all content, and takes full responsibility for the manuscript.

\end{document}